\documentclass{amsart}
\newtheoremstyle{bfnote}%
{}{}%
{\upshape}{}%
{\bfseries}{}%
{ }%
{\thmname{#1}\thmnumber{ #2}\thmnote{ (#3)}}
 
\usepackage{url}
\usepackage{amsmath, amssymb, amsthm, textcomp, rotating, float} 
\usepackage{hyperref}
\usepackage[capitalize,nameinlink]{cleveref}
\usepackage{mathrsfs}
 
\theoremstyle{bfnote}
\newtheorem{thm}{Theorem}[section]
\newtheorem{cor}[thm]{Corollary}
\newtheorem{lem}[thm]{Lemma}
\newtheorem{defi}[thm]{Definition}
\newtheorem{ex}[thm]{Example}
\newtheorem{prop}[thm]{Proposition}
\newtheorem{rmk}[thm]{Remark}

\newcommand{\Hom}{\text{Hom}}

\newcommand{\Der}{\text{Der}}
\renewcommand{\hat}{\widehat}
\newcommand{\id}{\text{id}}

\renewcommand{\phi}{\varphi}
\newcommand{\del}{\partial}

\renewcommand{\tilde}{\widetilde}

\newcommand{\Diff}{\text{Diff}}
\newcommand{\ddt}{\left.\frac{d}{dt}\right|_{t=0}}
\newcommand{\Def}{\text{def}}
\newcommand{\Hoch}{\text{Hoch}}
\newcommand{\GL}{\text{GL}}
\newcommand{\supp}{\text{supp}}
\newcommand{\chainmap}{\Phi}

\usepackage{xcolor,cancel}

\newcommand{\Gad}{\mathcal{G}_{\text{ad}}}

\title{Lie groupoid deformations and convolution algebras}
\author{Bjarne Kosmeijer}
\address{University of Amsterdam }
\email{b.a.kosmeijer@uva.nl}
\author{Hessel Posthuma}
\address{University of Amsterdam }
\email{h.b.posthuma@uva.nl}

\begin{document}

\maketitle

\begin{abstract}
We define a morphism from the deformation complex of a Lie groupoid to the Hochschild complex
of its convolution algebra, and show that it maps the class of a geometric deformation
to the algebraic class of the induced deformation in Hochschild cohomology. 
Applied to the adiabatic groupoid, we show that the van Est map to deformation cohomology of Lie algebroids is induced by taking the classical limit of a quantization map on the dual of the Lie algebroid.
\end{abstract}
\tableofcontents
\section*{Introduction}
The theory of Lie groupoids can be viewed as a blend of geometry and Lie theory, and plays an 
important role in several branches of mathematics. Lie groupoids can be viewed as the global 
integrations of geometrically defined Lie brackets and as such posses a natural deformation theory. 
In \cite{cms}, the authors wrote down an explicit cochain complex controlling such deformations. 
As expected, this ``deformation cohomology'' is isomorphic to the cohomology of the adjoint representation (up to homotopy as in \cite{ac}), but the point of \cite{cms} is that the latter involves the choice of a connection, whereas the deformation complex is intrinsically defined. Any deformation
of Lie groupoids defines a class in this deformation cohomology.

Lie groupoids also play an important role in noncommutative geometry where they give prime examples in the form of their associated convolution algebras. As is well-known, the deformation theory of an algebra is controlled by its Hochschild cohomology. Since a deformation of the underlying Lie groupoid induces a deformation of its convolution algebra, this strongly suggests a 
relationship between Hochschild cohomology and the adjoint representation.

The aim of the present article is to shed light on this issue by exhibiting an explicit morphism of 
cochain complexes between the deformation complex of a Lie groupoid and the Hochschild complex 
of its convolution algebra which relates the deformation classes associated to a deformation of the 
underlying Lie groupoid. We expect this morphism to be part of a larger picture computing the 
Hochschild cohomology in terms of higher powers of the adjoint representation.

A classical theme in the theory of Lie groupoids is its relation with the infinitesimal theory of Lie 
algebroids. For the deformation complex this is highlighted by the ``van Est'' morphism to the 
deformation complex of the Lie algebroid of Lie groupoid, a complex first considered in \cite{cm}.  

From the point of view of noncommutative geometry, the relation with the infinitesimal theory follows 
from ``quantization and the classical limit'': we show that our cochain morphism can be extended to the adiabatic groupoid 
of \cite{connes} interpolating between a Lie groupoid and its Lie algebroid. The van Est-map is 
then obtained by constructing a quantization map on the dual of the Lie algebroid using an exponential map on the Lie groupoid. The picture is completed by the relationship between 
the deformation cohomology a Lie algebroid  and the Poisson cohomology of its dual. 

This article is organized as follows: in \cref{prelim} we recall the basic set-up of Lie groupoids and their convolution algebras. In \cref{chainmap} we construct the morphism between the deformation complex of a Lie groupoid and the Hochschild complex of its convolution algebra, and explore some of its properties. Finally, \cref{quant} is devoted to obtaining the van Est map using the adiabatic groupoi together with a quantization.

\section{Preliminaries}\label{prelim}
\subsection{Densities along the fibers of a submersion}
\label{das}
Let $V$ be an n-dimensional vector space. A \textit{density} of $V$ is a map $a:\Lambda^n V\to\mathbb{R}$ such that for every invertible map $A\in\GL(n,\mathbb{R})$ it holds that $a(Av_1,...,Av_n)=|\det(A)|a(v_1,...,v_n)$.

More generally, from a vector bundle $E\to M$, one constructs a bundle of densities $\mathcal{D}_E$. Then if one has a vector bundle isomorphism $\Psi: E\to E$ covered by a diffeomorphism $\Phi: M\to M$, one obtains an action on the section of $\mathcal{D}_E$, defined by
\begin{equation*}
(\Psi^\ast a)_x(v_1,...,v_n)=a_{\Phi(x)}(\Psi v_1,...,\Psi v_n)
\end{equation*}
The case where $E=TM$ is of particular interest because the integral $\int_Ma$ is 
canonically defined for a compactly supported density of $TM$. In this case one obtains an action of a vector field $X\in\mathfrak{X}(M)$ on the densities on $TM$, namely:
\begin{equation}
\label{apvf}
Xa=\ddt (\Phi^t_X)^\ast a,
\end{equation}
where $\Phi^t_X$ denotes the flow of $X$. 

We will be mostly interested in densities along the fibers of a submersion. For this, let $f: M\to N$ be a submersion, and denote by $\mathcal{D}_f$ the bundle of densities of the vector bundle $\ker df$. In this case the fiber integral
\[
\int_f:\Gamma_c(M,\mathcal{D}_f)\to C^\infty_c(N)
\]
is canonically defined. A vector field $X$ acts on sections of $\mathcal{D}_f$ provided that the flow preserves the fibers of $f$. This is equivalent to there being a vector field $Y\in\mathfrak{X}(N)$ such that $df\circ X=Y\circ f$. In this case $X$ is called \textit{$f$-projectable}, and since $\Phi^t_X\circ f=f\circ\Phi^t_Y$ the flow of $X$ preserves the fibers of $f$ and in turn acts on $\ker df$, and we obtain an action of $X$ on $\Gamma(\mathcal{D}_f)$ by formula \eqref{apvf}.
We denote by $\Diff_f(M)$ the diffeomorphisms of $M$ that preserve the fibers of $f$, and by $\mathfrak{X}_f(M)$ the $f$-projectable vector fields of $M$.

In the following we shall consider $f$-projectable vector fields defined on only a single fiber of $f$ and let it act on densities to get a density on that one fiber. This is similar to the fact that the directional derivative $X(f)(p)$ of a function $f$ along a vector field $X$ in a point $p$ only depends
on $X(p)$. 
\begin{lem}
\label{nl}
Let $a\in\Gamma(\mathcal{D}_f)$, let $y\in N$, and let $X\in\mathfrak{X}_f(M)$ be an $f$-projective vector field. If $X$ vanishes along $f^{-1}(y)$, then $(Xa)_x=0$ for all $x\in f^{-1}(y)$.
\end{lem}
\begin{proof}
If $X$ vanishes along $f^{-1}(y)$ we have $\Phi^t_X(x)=x$ for all $t$ and all $x\in f^{-1}(y)$. In particular $d(\Phi^t_X)_x(v)=v$ for all $v\in\text{ker}(df)\subset T_xM$. This means that $((\Phi^t_X)^\ast a)_x=a_x$ and hence $(Xa)_x=0$.
\end{proof}
\begin{rmk}
The previous Lemma allows us to define $(Xa)_x$ for $x\in f^{-1}(y)$, $a\in\Gamma(\mathcal{D}_f)$ and $X\in\mathfrak{X}_f(M)|_{f^{-1}(y)}$. Indeed, we can choose $Y\in\mathfrak{X}_f(M)$ to be an extension of $X$ to a global vector field and define $(Xa)_x=(Ya)_x$. The previous Lemma is then used to show that this definition is independent of the choice of $Y$.
\end{rmk}
\subsection{The convolution algebra of a Lie groupoid}
Let $\mathcal{G}\rightrightarrows M$ be a Lie groupoid. For an introduction to the theory of Lie groupoids we refer to \cite{mm}. Here we denote source and target maps by $s,t:\mathcal{G}\to M$ and will think of arrows $g\in\mathcal{G}$ as pointing from right to left, so that the product $g_1g_2$ is defined whenever $s(g_1)=t(g_2)$. The Lie algebroid  $A(\mathcal{G})$ is defined as $A(\mathcal{G})=\ker(ds)|_M$. We will concern ourselves with the convolution algebra of $\mathcal{G}$. To define the convolution product, we need entities which can be integrated, and this is where densities come into play. To this end we look at densities along the source-fibers, where we note that there is a canonical isomorphism between $\ker ds$ and $t^\ast A(\mathcal{G})$ using right translations. In this way we can define the convolution product for two compactly supported densities $a_1, a_2\in\Gamma_c(\mathcal{D}_s)$ by
\begin{equation*}
(a_1\ast a_2)_g(v_1,...,v_n)=\int_{h\in s^{-1}(s(g))}(a_1)_{gh^{-1}}(v_1,...,v_n)(a_2)_h
\end{equation*}
In this notation $v_1,...,v_n\in A_{t(g)}=A_{t(gh^{-1})}$ so that the product in the integrand yields a well-defined compactly supported density along $s^{-1}(h)=s(g)$ that can be integrated. Colloquially this product will be written as:
\begin{equation*}
(a_1\ast a_2)(g)=\int_{g_1g_2=g}a_1(g_1)a_2(g_2)=\int_{h\in s^{-1}(s(g))}a_1(gh^{-1})a_2(h)
\end{equation*}
We define the convolution algebra $\mathcal{A}_\mathcal{G}$ of $\mathcal{G}$ to be $\mathcal{A}_\mathcal{G}=(\Gamma_c(\mathcal{D}_s),\ast)$. This definition of the convolution algebra differs slightly (but is isomorphic as a complex algebra) from the more usual one in e.g. \cite{connes} using $1/2$-densities along source {\em and} target fibers.
\subsection{The deformation complex of a Lie groupoid}
 Let $\mathcal{G}\rightrightarrows M$ be a Lie groupoid and write $\overline{m}$ for the map $\overline{m}(g,h)=gh^{-1}$. Note that $\overline{m}$ has as domain $\mathcal{G}\hspace*{1mm}^s\!\times^s\mathcal{G}:=\{(g_1,g_2)\in\mathcal{G}\times\mathcal{G},~s(g_1)=s(g_2)\}$. Furthermore we write $\mathcal{G}^{(k)}$ for the $k$'th nerve of $\mathcal{G}$:
\begin{equation*}
\mathcal{G}^{(k)}=\{(g_1,...,g_k)\in\mathcal{G}^k|s(g_i)=t(g_{i+1})\}
\end{equation*}
In \cite{cms} the deformation complex is defined as follows:
\begin{defi}
For $k\geq 1$ define 
$C^k_\Def(\mathcal{G})$ to be the set of smooth maps $c:\mathcal{G}^{(k)}\to T\mathcal{G}$ such that $c(g_1,...,g_k)\in T_{g_1}\mathcal{G}$ and such that there is a section $s_c$ of tthe vector bundle $t^*TM$ over $\mathcal{G}^{(k-1)}$ such that 
\[
ds(c(g_1,...,g_k))=s_c(g_2,...,g_k).
\]
The differential $\delta: C^k_\Def(\mathcal{G})\to C^{k+1}_\Def(\mathcal{G})$ is defined by setting:
\begin{align*}
(\delta c)(g_1,...,g_{k+1})=&-d\overline{m}(c(g_1g_2,g_3,...,g_{k+1}),c(g_2,...,g_{k+1}))\\
&+\sum_{i=2}^k(-1)^ic(g_1,...,g_ig_{i+1},...,g_{k+1})
+(-1)^{k+1}c(g_1,...,g_k).
\end{align*}
The {\em deformation complex} is defined by the graded vector space $C^\bullet_\Def(\mathcal{G}):=\bigoplus_{k\geq 1} C^k_\Def(\mathcal{G})$ equipped with the differential $\delta$, its cohomology is denoted $H^\bullet_\Def(\mathcal{G})$. 
\end{defi}
\begin{rmk}
\label{deg0}
It is possible, as in \cite{cms}, to extend the deformation complex in degree zero by putting $C^0_\Def(\mathcal{G})=\Gamma(M,A(\mathcal{G}))$ 
with differential defined for $\alpha\in\Gamma(M,A(\mathcal{G}))$ by
\begin{equation*}
(\delta \alpha)(g)=(dr_g)(\alpha(t(g))+(d(l_g\circ\iota))(\alpha(s(g))
\end{equation*}
We exclude these elements in degree $0$ because, as we will see, these element cannot correspond to Hochschild $0$-cochains.
\end{rmk}

\begin{rmk}
\label{mv}
It follows from the definition above that the closed elements in degree $1$ are exactly the multiplicative vector fields, c.f.\ \cite[\S 4.3]{cms}. These are vector fields $X\in\mathfrak{X}(\mathcal{G})$ that are $s$ and $t$-projectable to the same image in $\mathfrak{X}(M)$, satisfying the following equation:
\begin{equation*}
dm_{(g,h)}(X(g),X(h))=X(gh)
\end{equation*}
\end{rmk}

For certain purposes, most importantly applying the Van Est map, it often necessary to impose more strict relations on elements $c\in C^k_\Def(\mathcal{G})$ and their symbol $s_c$. To this end we also introduce the normalized deformation complex:

\begin{defi}
The \textit{normalized deformation complex} is the subcomplex $\hat{C}^\bullet_\Def(\mathcal{G})$ of $C^\bullet_\Def(\mathcal{G})$ consisting of those elements $c\in C^k_\Def(\mathcal{G})$ which satisfy
\begin{equation*}
c(1_x,g_2,...,g_k)=du(s_c(g_2,...,g_k))
\end{equation*}
and
\begin{equation*}
s_c(g_2,...,1_x,...,g_k)=0
\end{equation*}
where the unit is put in any of the $k-1$ slots.
\end{defi}
It is shown in \cite[Prop 11.8]{cms} that the inclusion of the normalized deformation complex into the whole deformation complex is a quasi-isomorphism.
\section{From deformation to Hochschild cohomology}\label{chainmap}
\subsection{The cochain map}
\label{cmh}
In this section we define a cochain map from the deformation complex of $\mathcal{G}$ to the Hochschild complex of the convolution algebra $\mathcal{A}_\mathcal{G}$. As a first hint for the existence of such a morphism, we make the following observation:
\begin{prop}\label{multvfderiv}
Let $\mathcal{G}\rightrightarrows M$ be a Lie groupoid. 
Multiplicative vector fields on $\mathcal{G}$  act as derivations on the convolution algebra.
\end{prop}
\begin{proof}
Recall the definition of a multiplicative vector field from \cref{mv}. Since a multiplicative vector field on $\mathcal{G}$ is by definition $s$-projectable to $M$, its action on an $s$-density is well-defined by 
the discussion in \cref{das}, c.f.\ equation \eqref{apvf}.

The key ingredient in the proof is the observation that the flow of a multiplicative vector field is a groupoid map, that is if $X\in\mathfrak{X}(\mathcal{G})$ is a multiplicative vector field then $\Phi^t_X(gh^{-1})=\Phi^t_X(g)\Phi^t_X(h)^{-1}$. A simple calculation then shows
\begin{align*}
X(a_1 \ast a_2)(g)&=\left.\frac{d}{dt}\right|_{t=0}(a_1\ast a_2)(\Phi^t_Xg)\\
&=\ddt\int_{h\in s^{-1}(s(\Phi^t_X g))}a_1((\Phi^t_X g)h^{-1})a_2(h)\\
&=\ddt\int_{h\in s^{-1}(s(g))}a_1((\Phi^t_X g)(\Phi^t_X h)^{-1})a_2(\Phi^t_Xh)\\
&=\ddt\int_{h\in s^{-1}(s(g))}a_1(\Phi^t_X(gh^{-1}))a_2(\Phi^t_Xh)\\
&=\ddt\int_{h\in s^{-1}(s(g))}a_1(\Phi^t_X(gh^{-1}))a_2(h)\\ &\hspace{2.5cm}+\ddt\int_{h\in s^{-1}(s(g))}a_1(gh^{-1})a_2(\Phi^t_Xh)\\
&=(Xa_1\ast a_2)(g)+(a_1\ast Xa_2)(g)
\end{align*}
which proves the proposition.
\end{proof}
In the following we write $C^\bullet_\Hoch(\mathcal{A}_\mathcal{G},\mathcal{A}_\mathcal{G})$ for the Hochschild complex of the convolution algebra $\mathcal{A}_\mathcal{G}$ with values in the bimodule $\mathcal{A}_\mathcal{G}$ with differential $\delta_{\rm Hoch}$.
We now describe the cochain map $C^\bullet_\Def(\mathcal{G})\to C^\bullet_\Hoch(\mathcal{A}_\mathcal{G},\mathcal{A}_\mathcal{G})$. 
\begin{defi}
\label{chainmap}
The map $\chainmap:C^\bullet_\Def(\mathcal{G})\to C^\bullet_\Hoch(\mathcal{A}_\mathcal{G},\mathcal{A}_\mathcal{G})$ is defined by 
\begin{equation*}
(\chainmap c)(a_1,...,a_k)(g)=\int_{g_1\cdots g_k=g}(c(-,g_2,...,g_k)a_1)(g_1)a_2(g_2)\cdots a_k(g_k),\qquad\mbox{for}~c\in C^k_\Def(\mathcal{G}).
\end{equation*}
This formula should be read as an inductive convolution (first over $g_1g_2=h_1$, then over $h_1g_3=h_2$, et cetera).
\end{defi}
\begin{rmk}
The formula for $\chainmap$ above is justified by \cref{nl}: $c\in C^k_\Def(\mathcal{G})$ is $s$-projectable and therefore the action of $c(-,g_2,...,g_k)$ on  $a\in \mathcal{A}_\mathcal{G}$ along $s^{-1}(t(g_2))$ is well-defined. In particular $(c(-,g_2,...,g_k)a)(g_1)$ is a well-defined density at $g_1$ for $(g_1,...,g_k)\in\mathcal{G}^{(k)}$.
\end{rmk}
Showing that $\chainmap$ is a chain map is done by a calculation similar to the one in \cref{multvfderiv}. In particular we need to deal with the term $\Phi^t_X(g)(\Phi^t_X(h))^{-1}$ for divisible $g$ and $h$ when $t$ goes to $0$. In the mulitplicative case this is precisely $\Phi^t_X(gh^{-1})$, but for general deformation elements we need a more general description.

For this we abbreviate the term in $\delta c$ involving $d\overline{m}$ by $\overline{m}c$, that is:
\begin{equation*}
(\overline{m}c)(g_1,...,g_{k+1})=d\overline{m}(c(g_1g_2,...,g_{k+1}),c(g_2,...,g_{k+1}))
\end{equation*}
We remark that this notation commutes with keeping the last all-but-two entries fixed, i.e.:
\begin{equation*}
\overline{m}(c(-,g_3,...,g_{k+1}))(g_1,g_2)=(\overline{m}c)(g_1,...,g_{k+1})
\end{equation*}
The key Lemma is then as follows:
\begin{lem}
Let $x\in M$, $X\in\mathfrak{X}_s(\mathcal{G})|_{s^{-1}(x)}$ and $a_1,a_2\in \mathcal{A}_\mathcal{G}$. Then for all $h\in s^{-1}(x)$ we have $\overline{m}X(-,h)\in\mathfrak{X}_s(\mathcal{G})|_{s^{-1}(t(h))}$ and for $g\in s^{-1}(x)$:
\begin{equation*}
X(a_1\ast a_2)(g)=(a_1\ast Xa_2)(g)+\int_{h\in s^{-1}(x)}((\overline{m}X(-,h))a_1)(gh^{-1})a_2(h)
\end{equation*}
\end{lem}
\begin{proof}
By definition we have:
\begin{equation*}
\overline{m}X(gh^{-1},h)=d\overline{m}(X(g),X(h))\in T_{gh^{-1}}\mathcal{G}
\end{equation*}
with $s$-projection
\begin{equation*}
ds(\overline{m}X(gh^{-1},h))=dt(X(h))
\end{equation*}
so indeed $\overline{m}X(-,h)\in\mathfrak{X}_s(\mathcal{G})|_{s^{-1}(t(h))}$.

Next we assume that $X$ is a globally defined $s$-projectable vector field (otherwise, we choose an extension at this point). Then we know that $\overline{m}X(-,h)$ is generated by the path $\Phi_t$ through $\Diff_s(\mathcal{G})$, which along $s^{-1}(t(h))$ looks like:
\begin{equation*}
\Phi_t(gh^{-1})=\Phi^t_X(g)(\Phi^t_X(h))^{-1}
\end{equation*}
so that we see that:
\begin{equation*}
\int_{h\in s^{-1}(x)}((\overline{m}X(-,h))a_1)(gh^{-1})a_2(h)=\ddt\int_{h\in s^{-1}(x)}a_1(\Phi^t_X(g)(\Phi^t_X(h))^{-1})a_2(h)
\end{equation*}
Using this we calculate $X(a_1\ast a_2)(g)$:
\begin{align*}
X(a_1\ast a_2)(g)=&\ddt (a_1\ast a_2)(\Phi^t_Xg)\\
=&\ddt\int_{h\in s^{-1}(s(\Phi^t_Xg))}a_1((\Phi^t_Xg)h^{-1})a_2(h)\\
=&\ddt\int_{h\in s^{-1}(x)}a_1((\Phi^t_Xg)(\Phi^t_Xh)^{-1})a_2(\Phi^t_X(h))\\
=&\ddt\int_{h\in s^{-1}(x)}a_1((\Phi^t_Xg)(\Phi^t_Xh)^{-1})a_2(h)+\ddt\int_{h\in s^{-1}(x)}a_1(gh^{-1})a_2(\Phi^t_Xh)\\
=&\int_{h\in s^{-1}(x)}((\overline{m}X(-,h))a_1)(gh^{-1})a_2(h)\\
&+(a_1\ast Xa_2)(g)
\end{align*}
which finishes the proof.
\end{proof}
\begin{prop}
The map $\chainmap: C^\bullet_\Def(\mathcal{G})\to C^\bullet_\Hoch(\mathcal{A}_\mathcal{G},\mathcal{A}_\mathcal{G})$ is a morphism of cochain complexes.
\end{prop}
\begin{proof}
This proof is essentially writing out all the parts of the Hochschild differential and applying some bookkeeping. We start with $c\in C^k_\Def(\mathcal{G})$ for $k\geq 1$, and write down the definition of the various parts of $\delta_\text{Hoch}(\chainmap c)$.
\begin{equation*}
(a_1\ast(\chainmap c)(a_2,...,a_{k+1}))(g)=\int_{g_1\cdots g_{k+1}=g}a_1(g_1)(c(-,g_3,...,g_{k+1})a_2)(g_2)a_3(g_3)\cdots a_{k+1}(g_{k+1})
\tag{$\star$}
\end{equation*}
\begin{equation*}
-(\chainmap c)(a_1\ast a_2,a_3,...,a_{k+1})(g)=-\int_{h\cdot g_3\cdots g_{k+1}=g}(c(-,g_3,...,g_{k+1})(a_1\ast a_2))(h)a_3(g_3)\cdots a_{k+1}(g_{k+1})\tag{$\star\star$}
\end{equation*}
\begin{multline*}
\sum_{i=2}^k(-1)^i(\chainmap c)(a_1,...,a_i\ast a_{i+1},...,a_{k+1})(g)=\\
\sum_{i=2}^{k}\int_{g_1\cdots g_{k+1}=g}(-1)^i (c(-,g_2,...,g_ig_{i+1},...,g_{k+1})a_1)(g_1)a_2(g_2)\cdots a_{k+1}(g_{k+1})
\end{multline*}
\begin{equation*}
(-1)^{k+1}((\chainmap c)(a_1,...,a_k)\ast a_{k+1})(g)=(-1)^{k+1}\int_{g_1\cdots g_{k+1}=g}(c(-,g_2,...,g_k)a_1)(g_1)a_2(g_2)\cdots a_{k+1}(g_{k+1})
\end{equation*}
The latter two terms we recognize from the differential of the deformation complex, while the first two terms can be rewritten to:
\begin{multline*}
(\star)+(\star\star)=\int_{hg_3\cdots g_{k+1}=g}\left((a_1\ast (c(-,g_3,...,g_{k+1}) a_2)-c(-,g_3,...,g_{k+1})(a_1\ast a_2)\right)(h)a_3(g_3)\cdots a_{k+1}(g_{k+1})
\end{multline*}
Then by the key Lemma we can rewrite this to
\begin{align*}
(\star)+(\star\star)=&-\int_{g_1\cdots g_{k+1}=g}((\overline{m}c)(-,g_2,...,g_{k+1})a_1)(g_1)a_2(g_2)\cdots a_{k+1}(g_{k+1})
\end{align*}
Putting this all together we conclude that:
\begin{align*}
(\delta_\Hoch(\chainmap c))(a_1,...,a_{k+1})(g)
&=(\chainmap(\delta c))(a_1,...,a_{k+1})(g)
\end{align*}
So we see that $\chainmap$ is indeed a chain-map.
\end{proof}
\subsection{Comparing deformation classes}
In this section we compare the deformation classes in $H^2_\Def(\mathcal{G})$ and $H^2_{\rm Hoch}(\mathcal{A}_\mathcal{G})$ coming from deformations of the Lie groupoid $\mathcal{G}$.
Recall from \cite[\S 5.2]{cms} that an $s$-constant deformation of $\mathcal{G}$ is a smooth family $\overline{m}_\epsilon: \mathcal{G}\hspace*{1mm}^s\!\times^s\mathcal{G}\to\mathcal{G}$ of division maps parameterized by $\epsilon$ in an open interval in $\mathbb{R}$ containing $0$, such that $\overline{m}_0=\overline{m}$. This induces a deformation cocycle $\beta\in C^2_\Hoch(\mathcal{A}_\mathcal{G},\mathcal{A}_\mathcal{G})$ by deforming the associative algebra that is the convolution algebra:
\begin{equation*}
\beta(a_1,a_2)(g)=\left.\frac{d}{d\epsilon}\right|_{\epsilon=0}\int_{h\in s^{-1}(s(g))}a_1(\overline{m}_\epsilon(g,h))a_2(h)
\end{equation*}
On the other hand the deformation also induces a deformation element $\xi\in C^2_\Def(\mathcal{G})$ set for $(g,g')\in\mathcal{G}^{(2)}$ by:
\begin{equation*}
\xi(g,g'):=\left.\frac{d}{d\epsilon}\right|_{\epsilon=0}\overline{m}_\epsilon(gg',g').
\end{equation*}
By \cite[Lemma 5.3]{cms}, this cochain is is closed: $\delta\xi=0$.
\begin{prop}
The chain map $\chainmap$ sends $\xi$ to $\beta$.
\end{prop}
\begin{proof}
This follows from the observation that if $s(h)=s(g)$, then
\begin{equation*}
\xi(gh^{-1},h)=\left.\frac{d}{d\epsilon}\right|_{\epsilon=0}\overline{m}_e(g,h)
\end{equation*}
With this we see that
\begin{equation*}
\beta(a_1,a_2)(g)=\int_{h\in s^{-1}(s(g))}(\xi(-,h)a_1)(gh^{-1})a_2(h)=\chainmap(\xi)(a_1,a_2)(g),
\end{equation*}
exactly as needed.
\end{proof}
\begin{rmk}
In \cite[Prop 5.12]{cms} a deformation cocycle $\xi\in C^2_\Def(\mathcal{G})$ is assigned to any deformation (in particular those who are not $s$-constant), whose cohomology class is canonical. Then $\chainmap(\xi)$ induces a Hochschild cohomology class of degree 2, which is not immediately linked to a deformation of the convolution product, since if the source map changes the underlying space of the convolution algebra also changes as it consists of densities along the $s$-fibers. Indeed, in \cite{cms} the authors need an auxillary choice of a vector field on the larger deformation space to define the cocycle. This choice of an auxillary vector field is precisely what is needed to compare the various convolution algebras when the source map varies, and in this way $\chainmap$ maps $[\xi]\in H^2_\Def(\mathcal{G})$  to the Hochschild class of the deformation of the convolution product thus defined. 
\end{rmk}

\subsection{Compatibility with the characteristic map to cyclic cochomology}
Denote by $(C^\bullet_{\rm diff}(\mathcal{G}),\delta)$ the cochain complex of inhomogeneous groupoid cochains given by $C^k_{\rm diff}(\mathcal{G}):=C^\infty(\mathcal{G}^{(k)})$ with differential
\begin{align*}
\delta\varphi(g_1,\ldots,g_{k+1})&=\varphi(g_2,\ldots,g_{k+1})+\sum_{i=1}^k(-1)^i\varphi(g_1,\ldots,g_ig_{i+1},\ldots, g_k)+(-1)^{k+1}\varphi(g_1,\ldots, g_k).
\end{align*}
We can turn this cochain complex into a DGA by introducing the product $\cup:C^k_{\rm diff}(\mathcal{G})\times C^l_{\rm diff}(\mathcal{G})\to C^{k+l}_{\rm diff}(\mathcal{G})$ given by
\[
(\varphi\cup\psi)(g_1,\ldots,g_{k+l}):=\varphi(g_1,\ldots,g_k)\psi(g_{l+1},\ldots, g_{k+l}).
\]
In \cite{cms} it is shown that by replacing $\varphi$ by a deformation cochain $c\in C^k_{\rm def}(\mathcal{G})$ in the above formula, $C^\bullet_{\rm def}(\mathcal{G})$ becomes a right module over $C^\bullet_{\rm diff}(\mathcal{G})$. On the other hand, in \cite{ppt}, the smooth groupoid cohomology was used to construct cyclic cocycles. In this section we shall that these two structures are compatible with each others under the cochain map $\Phi$ to Hochschild cohomology of \cref{cmh}. We start by re-writing the map to cyclic cohomology of \cite{ppt} in the following way.

First recall that the Hochschild cochain complex $C^\bullet(\mathcal{A}_{\mathcal{G}},\mathcal{A}_{\mathcal{G}})$ can be given a DGA structure by introducing the product $\cup:C^k(\mathcal{A}_{\mathcal{G}},\mathcal{A}_{\mathcal{G}})\times C^l(\mathcal{A}_{\mathcal{G}},\mathcal{A}_{\mathcal{G}})\to C^{k+l}(\mathcal{A}_{\mathcal{G}},\mathcal{A}_{\mathcal{G}})$
\[
(D\cup E)(a_1,\ldots,a_{k+l}):=D(a_1,\ldots,a_k)*E(a_{k+1},\ldots,a_{k+l}).
\]
Construct a map $\Phi_0:C^\bullet_{\rm diff} (\mathcal{G})\to C^\bullet(\mathcal{A}_\mathcal{G},\mathcal{A}_\mathcal{G})$ by
\begin{equation}
\label{dgm}
\Phi_0(\varphi)(a_1,\ldots,a_k)(g):=\int_{g_1\cdots g_k=g}\varphi(g_1,\ldots,g_k)a_1(g_1)\cdots a_k(g_k).
\end{equation}
\begin{lem}
The map $\Phi_0:(C^\bullet_{\rm diff} (\mathcal{G}),\delta,\cup)\to (C^\bullet(\mathcal{A}_\mathcal{G},\mathcal{A}_\mathcal{G}),\delta_{\rm Hoch},\cup)$ is a morphism of DGA's.
\end{lem}
\begin{proof}
This is a straightforward computation.
\end{proof}
With this Lemma we can also equip the Hochschild complex $C^\bullet(\mathcal{A}_\mathcal{G},\mathcal{A}_\mathcal{G})$ with a module structure over $C^\bullet_{\rm diff}(\mathcal{G})$ by using the cup-product on Hochschild cochains:
\[
(D\cup E)(a_1,\ldots, a_{k+l}):=D(a_1,\ldots,a_k)E(a_{k+1},\ldots,a_{k+l}).
\] 
Explicitly, this module structure is given by
\[
D\cdot\varphi:=D\cup\Phi_0(\varphi).
\]
We then have:
\begin{prop}
The cochain map $\chainmap: C^\bullet_\Def(\mathcal{G})\to C^\bullet_\Hoch(\mathcal{A}_\mathcal{G},\mathcal{A}_\mathcal{G})$ is a morphism of $C^\bullet_{\rm diff}(\mathcal{G})$-modules.
\end{prop}
\begin{proof}
Let us start with the following case: For $c\in C^k_\Def(\mathcal{G})$ and $f \in C^\infty(\mathcal{G})=C^1_{\rm diff}(\mathcal{G})$ we have
\begin{equation*}
\chainmap(c\cup f)(a_1,...,a_{k+1})=\chainmap(c)(a_1,...,a_k)\ast (f\cdot a_{k+1})
\end{equation*}
The claim follows by carefully writing out the definition
\begin{align*}
\chainmap(c\cup f)(a_1,...,a_{k+1})(g)&=\int_{g_1\cdots g_{k+1}=g}((c\cup f)(-,g_2,...,g_{k+1})a_1)(g_1)a_2(g_2)\cdots a_{k+1}(g_{k+1})\\
&=\int_{g_1\cdots g_{k+1}=g}(f(g_{k+1})c(-,g_2,...,g_k)a_1)(g_1)a_2(g_2)\cdots a_{k+1}(g_{k+1})\\
&=\int_{g_1\cdots g_{k+1}=g}(c(-,g_2,...,g_k)a_1)(g_1)a_2(g_2)\cdots a_k(g_k)\left(f(g_{k+1})a_{k+1}(g_{k+1})\right)\\
&=\int_{hg_{k+1}=g}\int_{g_1\cdots g_k=h}(c(-,g_2,...,g_k)a_1)(g_1)a_2(g_2)\cdots a_k(g_k)\left(f(g_{k+1})a_{k+1}(g_{k+1})\right)\\
&=\int_{hg_{k+1}=g}\chainmap(c)(a_1,...,a_k)(h)\left(f(g_{k+1})a_{k+1}(g_{k+1})\right)\\
&=\left(\chainmap(c)(a_1,...,a_k)\ast (f\cdot a_{k+1})\right)(g)
\end{align*}
Hence by induction we obtain
\begin{equation*}
\chainmap(c\cup (f_1\otimes\cdots\otimes f_l))(a_1,...,a_{k+l})=\chainmap(c)(a_1,...,a_k)\ast (f_1\cdot a_{k+1})\ast\cdots\ast (f_l\cdot a_{k+l})
\end{equation*}
Writing $f\ast a=\Phi_0(f)(a)$, we can rewrite this as
\[
\chainmap(c\cup (f_1\otimes\cdots\otimes f_l))=\chainmap (c)\cup\Phi_0(f_1)\cup\ldots\cup\Phi_0(f_l).
\]
From this the general statement of the proposition follows. 
\end{proof}

 Now, analogous to the action of vector fields on differential forms in geometry, the Hochschild cochains act on Hochschild chains by contraction:
\[
C^k(\mathcal{A}_\mathcal{G},\mathcal{A}_\mathcal{G})\times C_l(\mathcal{A}_\mathcal{G})\longrightarrow C_{l-k}(\mathcal{A}_\mathcal{G}),\qquad (D,a)\mapsto\iota_Da,
\]
given explicitly by
\[
\iota_D(a_0\otimes\ldots\otimes a_{k+l}):=a_0D(a_1,\ldots,a_k)\otimes a_{k+1}\otimes\ldots\otimes a_{k+l}.
\]
This action satisfies the properties
\begin{align*}
\iota_D\circ\iota_E&=\iota_{D\cup E}\\
[b,\iota_D]&=\iota_{\delta D}.
\end{align*}
The analogue of the Cartan formula for the ``Lie derivative''  $L_D:=B\circ \iota_D+\iota_D\circ B$ in noncommutative geometry also holds true on the level of Hochschild homology.

Next, recall from \cite{ppt} that when $\mathcal{G}$ is unimodular we can define a trace on the convolution algebra $\mathcal{A}_\mathcal{G}$ by 
\[
\tau(a):=\int_Ma\Omega,
\]
with on the right hand side $\Omega$ a $\mathcal{G}$-invariant section of the bundle $\mathcal{D}_{A^*}\otimes\mathcal{D}_{TM}$, and we use the duality $\mathcal{D}_A\times\mathcal{D}_{A^*}\to\mathbb{R}$ together with the isomorphism $\mathcal{D}_s|_M=\mathcal{D}_A$, to obtain a density on $M$ that can be integrated. With this trace (a degree $0$ cyclic cocycle), the cochain map 
\begin{equation}
\label{chain-diff}
\Psi_\tau:(C^\bullet_{\rm diff}(\mathcal{G}),\delta)\longrightarrow (C^\bullet(\mathcal{A}_\mathcal{G}),b_{\rm Hoch}),
\end{equation}
constructed in \cite{ppt} is simply given by $\Psi_\tau(c):=\iota_{\Phi_0(c)}\tau$.
\begin{cor}
Let $c\in C^k_{\rm def}(\mathcal{G})$ and $f\in C^l_{\rm diff}(\mathcal{G})$. Then the following
identity holds true:
\[
\iota_{\chainmap(c\cup f)}\tau=\iota_{\chainmap(c)}\Psi_\tau(f).
\]
\end{cor}
With this Corollary, we can construct new cyclic cocycles on the convolution algebra. First of all, if we start with a smooth 
groupoid cocycle $\varphi\in C^k_{\rm diff}(\mathcal{G})$, we obtain a Hochschild cocycle by applying $\Psi_\tau$ as in \eqref{chain-diff}. A small computation shows that this cocycle is closed under the $B$-differential, i.e., $B\Psi_\tau(\varphi)=0$, when $\varphi$ is cyclic:
\[
\varphi(g_1,\ldots,g_k)=(-1)^k\varphi((g_1\cdots g_k)^{-1},g_1,\ldots,g_{k-1})
\]
We can work out similar conditions for elements $c\in C^k_\Def(\mathcal{G})$, but they are more involved. For example, for $k=2$ we find
\begin{equation*}
(d\iota)(c(g,g^{-1}))=-c(g^{-1},g).
\end{equation*}

\begin{rmk}
It is proved in \cite[\S 9]{cms} that $H^\bullet_\Def(\mathcal{G})\cong H^\bullet(\mathcal{G},{\rm Ad})$, where ${\rm Ad}$ denotes the adjoint representation up to homotopy constructed in \cite{ac}. Taking into account the morphism \eqref{dgm}, this strongly suggests to relabel the morphism  of \cref{chainmap} as $\Phi_1$ and conjecture the existence of a map $\Phi_p:H^\bullet(\mathcal{G},{\rm Sym}^p({\rm Ad}))\to H^\bullet_{\rm Hoch}(\mathcal{A}_\mathcal{G},\mathcal{A}_\mathcal{G})$ extending the cases $p=0,1$ described in this paper. This would naturally fit with the infinitesimal theory (see also the next section) and the computation in \cite{blom} of the Hochschild cohomology of the universal enveloping algebra $\mathcal{U}(A)$ of the Lie algebroid $A$:
\[
H^\bullet_{\rm Hoch}(\mathcal{U}(A),\mathcal{U}(A))\cong \bigoplus_{p\geq 0} H^\bullet_{CE}(A,{\rm Sym}^pA).
\]
\end{rmk}
\subsection{The case $k=0$}
For the chain map between $C^\bullet_\Def (\mathcal{G})$ and $C^\bullet_\Hoch(\mathcal{A}_\mathcal{G},\mathcal{A}_\mathcal{G})$ we have just defined, a natural question is whether it can be extended to degree $k=0$, c.f. \cref{deg0}. For this, one must find a  a map  $\chainmap^0: \Gamma(A)\to \mathcal{A}_\mathcal{G}$ which extends the chain map $\chainmap$. This is only possible if $\chainmap(\delta(\alpha))\in\Der(\mathcal{A}_\mathcal{G})$ is an inner derivation for every $\alpha\in \Gamma(A)$.

Intuitively it is clear that is should not be always possible, since the derivation $\chainmap(\delta(\alpha))$ includes taking derivatives, while an inner derivation $\del_H(a)$ only includes integrations. The following example presents a concrete counterexample:

\begin{ex} Consider the pair groupoid $\mathbb{R}\times\mathbb{R}\rightrightarrows \mathbb{R}$. For this groupoid, a bundle of densities is trivialized by $|dx|$, so that every compactly supported density is of the form $f|dx|$ for a compactly supported smooth function $f$.  Furthermore, a section of the algebroid is simply a vector field $X\in\mathfrak{X}(\mathbb{R})$ and for this example we take $X=\frac{\partial}{\partial x}$. We have
\begin{equation*}
\delta(X)(x,y)=(X(x),X(y))
\end{equation*}
so that in this case $\delta(X)=\frac{\partial}{\partial x}+\frac{\partial}{\partial y}$. This vector field has flow
\begin{equation*}
\Phi^t_{\delta(X)}(x,y)=(x+t,y+t)
\end{equation*}
Next we consider $\chainmap\left(\delta\left(\frac{\partial}{\partial x}\right)\right)$, so we look at the action of $\frac{\partial}{\partial x}+\frac{\partial}{\partial y}$ on a density $f(x,y)|dx|\in \mathcal{A}_{\mathbb{R}\times\mathbb{R}}$. We see
\begin{equation*}
(\Phi^t_{\delta(X)})^\ast (f(x,y)|dx|)=f(x+t,y+t)|d(x+t)|=f(x+t,y+t)|dx|
\end{equation*}
So that:
\begin{equation*}
\chainmap(\delta(X))(f|dx|)=\left(\frac{\partial f}{\partial x}+\frac{\partial f}{\partial y}\right)|dx|
\end{equation*}
Now suppose that there is some $g|dx|\in \mathcal{A}_{\mathbb{R}\times\mathbb{R}}$, such that $\chainmap(\delta(X))=\partial_H(g|dx|)$. Then since always $\partial_H(g|dx|)(g|dx|)=0$, we see that:
\begin{equation*}
\frac{\partial g}{\partial x}+\frac{\partial g}{\partial y}=0
\end{equation*}
so that
\begin{equation*}
g(x+t,y+t)=g(x,y)
\end{equation*}
Since $g$ has to be compactly supported, the only possibility is that $g=0$, which is obviously not a solution to $\chainmap(\delta(X))=\partial_H(g|dx|)$. We conclude that $\chainmap(\delta(X))$ is not an inner derivation.
\end{ex}
In fact, using supports as an argument, we can deduce that $\chainmap(X)$ can never be an inner derivation for any $X\in\mathfrak{X}_s(\mathcal{G})$.
\begin{prop}
Let $D\in\Hom(\mathcal{A}_\mathcal{G},\mathcal{A}_\mathcal{G})$ be a non-zero Hochschild-1-cochain. If $D$ satisfies $\supp(Da)\subset\supp(a)$, then there is no $b\in \mathcal{A}_\mathcal{G}$ such that $D=[-,b]$.
\end{prop}
\begin{proof}
Suppose by contrary that there is a $b$ such that $D=[-,b]$. Let $g\in\mathcal{G}$ and let $a\in \mathcal{A}_\mathcal{G}$ be supported arbitrarily close to $g$. For $h\in t^{-1}(s(g))$ outside of the isotropy of $s(g)$ we obtain:
\begin{equation*}
(a\ast b)(gh)=\int_k a(gk^{-1})b(kh)\sim a(g)b(h)
\end{equation*}
where we use that $a$ is only non-zero close enough to $g$. For the other part of the commutator we have
\begin{equation*}
(b\ast a)(gh)=\int_k b(gk^{-1})a(kh)=0
\end{equation*}
Since there is no way to let $kh$ come arbitrarily close to $g$ since $h$ is not in the isotropy of $s(g)$.

Since $\supp (Da)\subset \supp(a)$ we see that $(a\ast b)(gh)$ also has to be supported arbitrarily close to $g$, so that $b$ is identically zero outside of the isotropy of $\mathcal{G}$.

If we look at $h$ an isotropy element of $\mathcal{G}$ we see that the second term acts like $b(ghg^{-1})a(g)$, so that we see that $b$ is invariant under conjugation. However, if $b$ is invariant under conjugation we conclude that $b\in Z(\mathcal{A}_\mathcal{G})$, which is in contradiction to the fact that $D$ is non-zero. We conclude that there is no $b$ that solves $D=[-,b]$.
\end{proof}
\begin{rmk}
It is possible to define the map $\Phi^0$ if one allows for distributions to be cochains of degree $0$, that is if one defines $C^0_\text{Hoch}(\mathcal{A}_\mathcal{G},\mathcal{A}_\mathcal{G}):=\Gamma^{-\infty}_c(\mathcal{D}_s)$.
\end{rmk}
\begin{cor} If $X\in C^1_\Def(\mathcal{G})$ is non-zero, then $\chainmap(X)$ can never be an inner derivation.
\end{cor}
\begin{proof}
This follows from the previous proposition by the observation that $\chainmap(X)$ is local since it involves taking derivatives and the fact that $\chainmap$ is easily observed to be injective.
\end{proof}
\subsection{Examples}
In this section we discuss how the chain map $\chainmap$ links the deformation cohomology of $\mathcal{G}$ and the Hochschild cohomology of $\mathcal{A}_\mathcal{G}$ in certain examples.

\begin{ex}[Trivial groupoid]
We consider the trivial groupoid $\mathcal{G}=M\rightrightarrows M$. On the density side we simply have $(\mathcal{A}_\mathcal{G},\ast)=(C^\infty_c(M),\cdot)$, with $H^\bullet_\Hoch(\mathcal{A}_\mathcal{G},\mathcal{A}_\mathcal{G})=\Lambda^\bullet\mathfrak{X}(M)$. At the side of the deformation complex we note that the $k$-nerve of the trivial groupoid is $M$ for every $k$ and $s$-projectability is a void property, so that for $k>0$ we have $C^k_\Def(\mathcal{G})=\mathfrak{X}(M)$, with differential alternating between the identity and the zero map:
\begin{equation*}
C^\bullet_\Def(\mathcal{G})=\left[ 0\to \mathfrak{X}(M)\xrightarrow{0}\mathfrak{X}(M)\xrightarrow{\rm id}\mathfrak{X}(M)\to\cdots\right]
\end{equation*}
So the deformation cohomology equals:
\begin{equation*}
H^k_\Def(\mathcal{G})\cong\left\{\begin{matrix}
\mathfrak{X}(M) &\text{ if } k=1\\
0 & \text{ else}
\end{matrix}\right.
\end{equation*}
The chain map $\chainmap: C^\bullet_\Def(\mathcal{G})\to C^\bullet_\Hoch(\mathcal{A}_\mathcal{G},\mathcal{A}_\mathcal{G})$ simply becomes:
\begin{equation*}
\chainmap(X)(f_1,...,f_k)=(Xf_1)\cdot f_2\cdots f_k
\end{equation*}
and we simply see that:
\begin{equation*}
H^k(\chainmap)=\left\{\begin{matrix}
\id &\text{ if }k=1\\
0 & \text{ else}
\end{matrix}\right.
\end{equation*}
We should also remark for this example that using the classical Hochschild--Kostant--Rosenberg theorem, we see that taking exterior powers of deformation elements we retrieve the whole Hochschild cohomology of $C^\infty_c(M)$.
\end{ex}
\begin{ex}[\'Etale groupoids]
In the case of an \'Etale groupoid $\mathcal{G}\rightrightarrows M$, we have $\mathcal{A}_\mathcal{G}=C^\infty_c(\mathcal{G})$, since the distribution $\ker(ds)$ is the trivial distribution. The convolution product in this case is commonly written as
\begin{equation*}
(f_1\ast f_2)(g)=\sum_{g_1g_2=g}f_1(g_1)f_2(g_2)
\end{equation*}
In this case the action of vector fields on densities is just the normal action of vector fields on functions, and the map $\chainmap$ reduces to
\begin{equation*}
\chainmap(c)(f_1,...,f_k)(g)=\sum_{g_1\cdots g_k=g}(c(g_1,...,g_k)f_1)\cdot f_2(g_2)\cdots f_k(g_k)
\end{equation*}
Since the source map of $\mathcal{G}$ is a local diffeomorphism, we see that there is a 1-1 correspondence between deformation elements $c\in C^k_\Def(\mathcal{G})$ and their symbols $s_c\in\Gamma(t^\ast TM\to\mathcal{G}^{(k-1)})$ since we have
\begin{equation*}
c(g_1,...,g_k)=(ds_{g_1})^{-1}(s_c(g_2,...,g_k))
\end{equation*}
In fact, the correspondence establishes an isomorphism between $C^\bullet_\Def(\mathcal{G})$ and $C^\bullet(\mathcal{G},TM)[-1]$ where we see $TM$ as a representation of $\mathcal{G}$ where $g$ acts $T_{s(g)}M\to T_{t(g)}M$ as
\begin{equation*}
g\cdot v=dt_g((ds_g)^{-1})(v))
\end{equation*}
The shift by 1 we see here also serves as a justification of why the case $k=0$ is a tricky thing (although for \'Etale groupoids of course we have $C^0_\Def(\mathcal{G})=0$).

In the case that we have a proper \'Etale groupoid (over a connected base $M$) we can calculate the cohomologies in both sides of the equation. On the side of the deformation complex we use \cite[Thm 6.1]{cms} to obtain:
\begin{align*}
H^0_\Def(\mathcal{G})&\cong\{0\}\\
H^1_\Def(\mathcal{G})&\cong\mathfrak{X}(M)_\text{inv}\\
H^k_\Def(\mathcal{G})&\cong\{0\}\hspace*{1cm}(k\geq 2)
\end{align*}
For the Hochschild cohomology of the convolution algebra we refer to \cite[Thm 3.11]{nppt} to obtain
\begin{align*}
H^k(\mathcal{A}_\mathcal{G},\mathcal{A}_\mathcal{G})\cong \bigoplus_{\mathcal{O}\in\text{Sec}(\mathcal{G})}\Gamma_\text{inv}(\Lambda^{k-\text{codim}(\mathcal{O})}T\mathcal{O})
\end{align*}
where the sum is over the sectors $\mathcal{O}$ of $\mathcal{G}$. The action of the chain map $\Phi$ on the cohomology of degree $1$ is the inclusion of $\mathfrak{X}(M)_\text{inv}$ into this sum as the term for the sector $\mathcal{O}=M$.
\end{ex}
\section{Deformation quantization and the van Est map}\label{quant}
\subsection{The adiabatic groupoid}\label{adiabatic}

In the theory of deformation quantizations and applications thereof, there is an inherent place for replacing a groupoid with its adiabatic groupoid, as first described in \cite{connes}. In view of our discussion of the deformation complex, we describe it using the division map:
\begin{defi}
Let $\mathcal{G}\rightrightarrows M$ be a Lie groupoid, with Lie algebroid $A\xrightarrow{\pi} M$. We define the \textit{adiabatic groupoid} $\Gad\to M\times\mathbb{R}$ by:
\begin{equation*}
\Gad=A\times\{0\}\sqcup \mathcal{G}\times\mathbb{R}^\ast
\end{equation*}
The source and target are defined by:
\begin{equation*}
s(v,0)=(\pi(v),0)
\end{equation*}
\begin{equation*}
s(g,\tau)=(s(g),\tau)
\end{equation*}
\begin{equation*}
t(v,0)=(\pi(v),0)
\end{equation*}
\begin{equation*}
t(g,\tau)=(t(g),\tau)
\end{equation*}
Then we define the inversion map by
\begin{equation*}
\iota(v,0)=(-v,0)
\end{equation*}
\begin{equation*}
\iota(g,\tau)=(\iota(g),\tau)
\end{equation*}
Lastly, to define we division map, we note that pairs of divisible arrows come in 2 shapes, namely pairs $(v,0)$ and $(w,0)$ with $\pi(v)=\pi(w)$, and pairs $(g,\tau)$ and $(h,\tau)$ where $g$ and $h$ are divisible. We then define the division map by:
\begin{equation*}
\overline{m}((v,0),(w,0))=(v-w,0)
\end{equation*}
\begin{equation*}
\overline{m}((g,\tau),(h,\tau))=(\overline{m}(g,h),\tau)
\end{equation*}
\end{defi}
This is just the set-theoretical description, but the remarkable feature is that the adiabatic groupoid can be given a smooth
structure. Here we briefly recall this smooth structure and show how to extend normalized deformation elements to deformation elements of the adiabatic groupoid. Both will be done in the context of the procedure known as the {\em deformation to the normal cone}.
\subsubsection{Deformation to the normal cone}
The part of the discussion below concerning the smooth structure and the smooth maps on the deformation to the normal cone is after \cite[\S 4]{higson} and \cite[\S 1.1]{ds}.
\begin{defi}\label{nms}
Let $S\hookrightarrow M$ be a submanifold with normal bundle $N\to S$. The \textit{deformation to the normal cone} $N(M,S)$ is the manifold defined by:
\begin{equation*}
N(M,S)=N\times\{0\}\sqcup M\times\mathbb{R}^\ast
\end{equation*}
\end{defi}
The deformation to the normal cone can be given a topology and smooth structure in two ways. Either it is characterized by the fact that the following two types of maps
\begin{itemize}
\item The map $N(M,S)\to M\times\mathbb{R}$ that sends $(x,\tau)$ for $\tau\neq  0$ to $(x,\tau)$ and sends $(v,0)$ with $v\in N_x$ to $(x,0)$.
\item For every $f\in C^\infty(M)$ such that $f|_S=0$, the map $\delta f: N(M,S)\to\mathbb{R}$ defined by
\begin{align*}
(\delta f)(x,\tau)&=\frac{f(x)}{\tau}\,\,\,(x\in M,\,\tau\neq 0),\\
(\delta f)(v,0)&=d_nf(v)\,\,\,(v\in N)
\end{align*}
\end{itemize}
are smooth. Here by $d_nf$ we mean the smooth map on $N$ that for $v\in TM|_S$ sends $[v]$ to $df(v)$ and which is well-defined since $f|_S=0$.

Equivalently, one uses an exponential map, that is a map $\theta: U\to M$ from an open neighbourhood $U\subset N$ of the zero-section, with the property that for all $p\in S$ and $v\in N_p$ it holds that
\begin{equation*}
\theta(0_p)=p,\qquad 
\left.\frac{d}{d\tau}\right|_{\tau=0}\theta(\tau v)=v ~\text{ mod }T_pS
\end{equation*}
The smooth structure on $N(M,S)$ can then also be characterized by the fact that the maps
\begin{equation*}
i_1: M\times\mathbb{R}^\ast\to N(M,S):\,\,(x,\tau)\mapsto (x,\tau)
\end{equation*}
\begin{equation*}
i_2: U'=\{(v,\tau)\in N\times\mathbb{R}: \tau v\in U\}\to N(M,S):\,\,\begin{matrix}(v,\tau)&\mapsto& (\theta(\tau v),\tau)\\
(v,0)&\mapsto& (v,0)\end{matrix}
\end{equation*}
are open smooth embeddings.

Important in considering deformations to normal cones is the action of $\mathbb{R}^\ast$ on $N(M,S)$, which is given by:
\begin{align*}
\lambda\cdot (x,\tau)&=(x,\lambda\tau),\\
\lambda\cdot (v,0)&=\left(\frac{v}{\lambda},0\right)
\end{align*}
where $\lambda,\tau\in\mathbb{R}^\ast$, $x\in M$ and $v\in N$.

We will describe how to extend a vector field on $M$, that is parallel to $S$, to a vector field on $N(M,S)$ that is invariant under the $\mathbb{R}^\ast$-action.

This will be done by writing down a vector field on the normal bundle and combining it with a vector field over $M\times\mathbb{R}^\ast$ to a discrete vector field on $N(M,S)$, and using an explicit description of the smooth functions on $N(M,S)$ to show that this is in fact  a {\em smooth} vector field.
\begin{defi} \cite{higson}
Let $X$ be a set and $\mathcal{F}=\{f_\alpha: X\to V_\alpha\}$ be a family of functions from $X$ into smooth manifolds. We say that a function $f: X\to\mathbb{R}$ is \textit{smoothly composed from the family $\mathcal{F}$} if there is a finite collection $(f_{\alpha_1},...,f_{\alpha_n})\subset\mathcal{F}$ and a smooth map $h: V_{\alpha_1}\times\cdots V_{\alpha_n}\to\mathbb{R}$ such that
\begin{equation*}
f(x)=h(f_{\alpha_1}(x),...,f_{\alpha_n}(x))
\end{equation*}
\end{defi}
The smooth structure of $N(M,S)$ then means that all smooth functions on $N(M,S)$ are smoothly composed of type of functions as described after \cref{nms}. If we then apply Taylors theorem we conclude the following.
\begin{lem}\label{smoothvfnms}
A discrete vector field $X$ on $N(M,S)$ is smooth if and only if for every $f\in C^\infty(M)$ with $f|_s=0$ and every $g\in C^\infty(M\times\mathbb{R})$ the maps $\delta f$ and $\tilde{g}\in C^\infty(N(M,S))$ defined by:
\begin{equation*}
\begin{matrix}
(\delta f)(x,\tau)=\frac{f(x)}{\tau}&\,\,\,(\tau\neq 0)\\
(\delta f)(v,0)=d_nf(v)&\,\,\,(v\in N)\\
\\
\tilde{g}(x,\tau)=g(x,\tau)&\,\,\,(\tau\neq 0)\\
\tilde{g}(v,0)=g(x,0)&\,\,\,(v\in N_x)
\end{matrix}
\end{equation*}
satisfy that $X(\delta f),X(\tilde{g})\in C^\infty(N(M,S))$.
\end{lem}
We start with writing down the vector field over $N$. This is the {\em linearization}, as also in \cite[\S 4.1]{az}, that we describe in detail below:
\begin{prop}\label{vectorfieldnormalbundle}
Let $S\hookrightarrow M$ be a submanifold with normal bundle $\pi:N\to S$ and $X\in\mathfrak{X}(M)$ a vector field that is parallel to $S$. Then:
\begin{itemize}
\item[\textbf{a)}] The map that sends a smooth function $f\in C^\infty(M)$ satisfying $f|_S=0$ to the map $d_nf\in C^\infty_\text{lin} (N)$ is a surjection onto $C^\infty_\text{lin}(N)$
\item[\textbf{b)}] If $f\in C^\infty(M)$ satisfies that $f|_S=0$ and $d_nf=0$, then $Xf$ satisfies that $d_n(Xf)=0$.
\item[\textbf{c)}] The maps $(X_N)_{\text{lin}}: C^\infty_\text{lin}(N)\to C^\infty(N)$ and $(X_N)_{\text{cst}}: C^\infty(S)\to C^\infty(N)$ defined by
\begin{equation*}
(X_N)_\text{lin}(d_nf)=d_n(Xf)
\end{equation*}
\begin{equation*}
(X_N)_\text{cst}(g)=X|_S (g)\circ \pi
\end{equation*}
define a smooth vector field $X_N\in\mathfrak{X}(N)$.
\end{itemize}
\end{prop}
\begin{proof} Working down the list:
\begin{itemize}
\item[\textbf{a)}] By using a partition of unity this reduces to the local case $M=\mathbb{R}^m\times\mathbb{R}^n$ with $S=\mathbb{R}^m\times\{0\}$. In this local case there is a canonical diffeomorphism between $M$ and $N$ and pushing a linear map on $N$ through this canonical diffeomorphism yields a smooth map on $M$ which normal derivative equals the linear map on $N$ we started with.
\item[\textbf{b)}] This is again a computation in the local case $M=\mathbb{R}^m\times\mathbb{R}^n$ with $S=\mathbb{R}^m\times\{0\}$. Write
\begin{equation*}
X=\sum_{i=1}^m \alpha_i(x,y)\frac{\partial}{\partial x_i}+\sum_{j=1}^n\beta_j(x,y)\frac{\partial}{\partial y_j}
\end{equation*}
The fact that $X$ is parallel to $S$ means that $\beta_j(x,0)=0$ for all $j=1,...,n$. The fact that $d_nf=0$ is equivalent to the fact $\frac{\partial f}{\partial y_j}(x,0)=0$ for all $j=1,...,n$. Then we have
\begin{equation*}
Xf=\sum_{i=1}^m \alpha_i\frac{\partial f}{\partial x_i}+\sum_{j=1}^n \beta_j\frac{\partial f}{\partial y_j}
\end{equation*}
So that for $k=1,...,n$ we have
\begin{equation*}
\frac{\partial (Xf)}{\partial y_k}=\sum_{i=1}^m\frac{\partial \alpha_i}{\partial y_k}\frac{\partial f}{\partial x_i}+\sum_{i=1}^m\alpha_i\frac{\partial^2 f}{\partial y_k\partial x_i}+\sum_{j=1}^n\frac{\partial\beta_j}{\partial y_k}\frac{\partial f}{\partial y_j}+\sum_{j=1}^n\beta_j\frac{\partial^2 f}{\partial y_k\partial y_j}
\end{equation*}
Then since respectively $\frac{\partial f}{\partial x_i}(x,0)=0$ (since $f(x,0)=0$), $\frac{\partial^2 f}{\partial y_k\partial x_i}(x,0)=\left(\frac{\partial}{\partial x_i}\frac{\partial f}{\partial y_k}\right)(x,0)=0$ (since $\frac{\partial f}{\partial y_k}(x,0)=0$), $\frac{\partial f}{\partial y_j}(x,0)=0$ (by assumption) and $\beta_j(x,0)=0$ (by assumption), we see that
\begin{equation*}
\frac{\partial(Xf)}{\partial y_k}(x,0)=0
\end{equation*}
which implies that $d_n(Xf)=0$.
\item[\textbf{c)}] First note that (by restriction) a smooth vector field $Y\in\mathfrak{X}(E)$ on a vector bundle $\pi: E\to M$ is the same as a pair of maps $Y_\text{lin}: C^\infty_\text{lin}(E)\to C^\infty(E)$ and $Y_\text{cst}: C^\infty(M)\to C^\infty(E)$ such that for all $f,g\in C^\infty(M)$ and $h\in C^\infty_\text{lin}$ it holds that
\begin{equation*}
Y_\text{cst}(fg)=(f\circ\pi)\cdot Y_\text{cst}(g)+(g\circ\pi)\cdot Y_\text{cst}(f)
\end{equation*}
\begin{equation*}
Y_\text{lin}((f\circ\pi)\cdot h)=(f\circ\pi)\cdot Y_\text{lin}(h)+h\cdot Y_\text{cst}(f)
\end{equation*}
We show that these properties hold for the maps $(X_N)_\text{cst}$ and $(X_N)_\text{lin}$.

First we note that $(X_N)_\text{lin}$ is well-defined by parts a) and b). To show that they define a smooth vector field we check for $f,g\in C^\infty(S)$
\begin{align*}
(X_N)_\text{cst}(fg)&=(X|_S(fg))\circ\pi=(f\cdot X|_S(g)+g\cdot X|_S(f))\circ\pi\\
&=(f\circ\pi)\cdot (X|_S(g)\circ\pi)+(g\circ\pi)\cdot (X|_S(f)\circ\pi)\\
&=(f\circ\pi)X_\text{cst}(g)+(g\circ\pi)X_\text{cst}(f)
\end{align*}
Secondly let $f\in C^\infty(S)$ and $h\in C^\infty_\text{lin}(N)$ given by $h=d_ng$ with $g\in C^\infty(M)$ such that $g|_S=0$. Then first we need to find $g'\in C^\infty(M)$ with $g'|_S=0$ such that $fh=d_n(g')$. This can be done by choosing an extension of $f$ which is `constant in the normal direction', which is only well-defined locally or if we choose an exponential map.

We resort to the local case $M=\mathbb{R}^m\times\mathbb{R}^n$ with $S=\mathbb{R}^m\times\{0\}$. Then the map $g'(x,y)=f(x)g(x,y)$ clearly satisfies that $d_ng'=fh$. Then writing $X$ in coordinates as
\begin{equation*}
X=\sum_{i=1}^m\alpha_i\frac{\partial}{\partial x_i}+\sum_{j=1}^n\beta_j\frac{\partial}{\partial y_j}
\end{equation*}
we have
\begin{equation*}
(Xg')(x,y)=\sum_{i=1}^m\alpha_i(x,y)\frac{\partial f}{\partial x_i}(x)g(x,y)+f(x)(Xg)(x,y)
\end{equation*}
so that we see
\begin{equation*}
\frac{\partial (Xg')}{\partial y_k}(x,0)=\sum_{i=1}^m\alpha_i(x,0)\frac{\partial f}{\partial x_i}(x)\frac{\partial g}{\partial y_k}(x,0)+\sum_{i=1}^m\frac{\partial \alpha_i}{\partial y_k}(x,0)\frac{\partial f}{\partial x_i}(x)g(x,0)+f(x)\frac{\partial (Xg)}{\partial y_k}(x,0)
\end{equation*}
Then $g(x,0)=0$ so that the middle term vanishes. Then recognizing terms we obtain
\begin{equation*}
d(Xg')_{(x,0)}\left(\frac{\partial}{\partial y_k}\right)=X|_S(f)(x)\cdot(dg)_x\left(\frac{\partial}{\partial y_k}\right)+f(x)\cdot d(Xg)_{(x,0)}\left(\frac{\partial}{\partial y_k}\right)
\end{equation*}
so that globalizing we have
\begin{align*}
(X_N)_\text{lin}((f\circ \pi)\cdot d_ng)&=(X_N)_\text{lin}(d_ng')\\
&=d_n(Xg')\\
&=(X|_S(f)\circ\pi)d_ng+(f\circ\pi)d(Xg)\\
&=(X_N)_\text{cst}(f)d_ng+(f\circ\pi)(X_N)_\text{lin}(d_ng)
\end{align*}
So we see that we obtain a smooth vector field $X_N\in\mathfrak{X}(N)$.
\end{itemize}
This completes the proof.
\end{proof}
We are now ready to define the $\mathbb{R}^\ast$-invariant extension of the vector field $X$.
\begin{prop}\label{vfonnms}
Let $S\hookrightarrow M$ be a submanifold with normal bundle $N\to S$. Let $X\in\mathfrak{X}(M)$ be a vector field that is parallel to $S$. Then the discrete vector field $X_\text{inv}$ on $N(M,S)$ defined by
\begin{align*}
X_\text{inv}(x,\tau)&=X(x),\quad(\tau\neq 0)\\
X_\text{inv}|_{N\times\{0\}}&=X_N
\end{align*}
is a smooth vector field $X_\text{inv}\in\mathfrak{X}(M,S)$ which is the unique vector field on $N(M,S)$ which equals $X$ on $M\times\mathbb{R}^\ast$ and the unique $\mathbb{R}^\ast$-invariant vector field on $N(M,S)$ which equals $X$ along $M\times\{1\}$.
\end{prop}
\begin{proof}
The invariance and uniqueness is clear assuming that $X_\text{inv}$ is smooth. To show that it is smooth, by \cref{smoothvfnms} the only thing we have to check is that $X_\text{inv}(\delta f)$ and $X_\text{inv}(\tilde{g})$ are smooth for $f\in C^\infty(M)$ with $f|_S=0$ and $g\in C^\infty(M\times\mathbb{R})$. The definition of $X_N$ makes sure that the result is
\begin{equation*}
X_\text{inv}(\delta f)=\delta(Xf)
\end{equation*}
\begin{equation*}
X_\text{inv}(\tilde{g})=\tilde{Xg}
\end{equation*}
where in the second equation $X$ acts on $C^\infty(M\times\mathbb{R})$ as the vector field $X(x,\tau)=X(x)$ on $M\times\mathbb{R}$. By definition $\delta(Xf)$ and $\tilde{Xg}$ are smooth and so the result follows.
\end{proof}
\subsubsection{The adiabatic groupoid as a deformation to the normal cone}
We can now apply this to the case $M\hookrightarrow\mathcal{G}$ with normal bundle $A=\ker ds|_M$. The fact that the source, target and division maps are smooth, follows from the fact that away from $\tau=0$ they are just the respective maps of the original groupoid, while along $\tau=0$ they are the normal derivatives of the respective maps. A general principle of deformations to normal cones then means they are smooth. We note that an exponential map can be obtained by choosing a connection on $A$, see \cite{nwx} and \cite{landsmanboek}.

Next we want to describe the nerve of the adiabatic groupoid. As a set it equals $(\Gad)^{(k)}=\mathcal{G}^{(k)}\times\mathbb{R}^\ast\sqcup A^{\oplus k}\times\{0\}$. From the view point of trying to define vector fields on the nerve of the adiabatic groupoid, this set-theoretic description leads to searching for a connection between $A^{\oplus k}$ and the normal bundle of $M$ inside $\mathcal{G}^{(k)}$ as the diagonal of units.

\begin{lem}\label{normalbundenerve}
Let $\mathcal{G}\rightrightarrows M$ be a Lie groupoid with $\Delta: M\to\mathcal{G}^{(k)}$ the diagonal inclusion via the units. The vector bundle map $\nu: A^{\oplus k}\to\Delta^\ast T\mathcal{G}^{(k)}$ given by
\begin{equation*}
\nu(v_1,...,v_k)=(v_1+\sum_{i=2}^k du(dt(v_i)),v_2+\sum_{i=3}^k du(dt(v_i)),...,v_{k-1}+du(dt(v_k)),v_k)
\end{equation*}
induces an isomorphism between $A^{\oplus k}$ and the normal bundle of $M$ inside $\mathcal{G}^{(k)}$.
\end{lem}
\begin{proof}
First one checks that $\nu$ indeed maps into the tangent space of $\mathcal{G}^{(k)}\subset\mathcal{G}^{\times k}$, which is a simple calculation. Next to show that it induces an isomorphism to the normal bundle to $\Delta$, we first use the decomposition $T_{1_x}M=A_x\oplus T_xM$ to see that if $\nu(v_1,...,v_k)\in T_xM\subset T_{\Delta(x)}\mathcal{G}^{(k)}$ then $(v_1,...,v_k)=0$, so that the map into the normal bundle is injective. A simple case of dimension counting then implies that it the induced map is an isomorphism.
\end{proof}
\begin{cor}
There is a natural isomorphism between $N(\mathcal{G}^{(k)},M)$ and $\Gad^{(k)}$ which away from $\tau=0$ links $((g_1,...,g_k),\tau)$ and $((g_1,\tau),...,(g_k,\tau))$.
\end{cor}
\subsubsection{Haar systems on the adiabatic groupoid}
We intend to link deformation quantizations of the Poisson manifold $A^\ast$ with the Van Est map $\mathcal{V}:\tilde{C}^\bullet_\Def(\mathcal{G})\to C^\bullet_\Def(A)$. To make the syntax line up, we need to explicitely write down isomorphisms between smooth functions on $\mathcal{G}$ and elements of the convolution algebra. This is done via Haar systems, which we will describe here in terms of densities.

\begin{defi}A Haar system on a groupoid $\mathcal{G}\rightrightarrows M$ is a collection $\lambda=\{\lambda_x\}_{x\in M}$ of positive sections $\lambda_x\in\Gamma(\mathcal{D}_s|_{s^{-1}(x)})$ that are invariant under right translations $R_g: s^{-1}(t(g))\to s^{-1}(s(g))$ and such that for every compactly supported function $f\in C^\infty_c(\mathcal{G})$ the map $\lambda(f): M\to\mathbb{R}$ given by
\begin{equation*}
\lambda(f)(x)=\int_{s^{-1}(x)}f(g)\lambda_x(g)
\end{equation*}
is smooth.
\end{defi}
We know that every Lie groupoid admits a Haar system (\cite[Prop 3.4]{landsman}) and if we have a Haar system $\lambda$ on a Lie groupoid $\mathcal{G}\rightrightarrows M$ with $s$-fibers of dimension $d$, we can (\cite[p.19]{landsman}) induce a Haar system $\hat{\lambda}$ on $\Gad$ given by
\begin{equation*}
\hat{\lambda}(g,\tau)=|\tau|^d\lambda(g)
\end{equation*}
\begin{equation*}
\hat{\lambda}(v,0)=\lambda(\pi(v))
\end{equation*}
Here $\pi: A\to M$ is the projection and we take the canonical isomorphism $\ker(d\pi)\cong\pi^\ast(A)$ as a given.

Note that in particular we obtain a Haar system on the vector bundle $A\to M$, seen as a groupoid in the canonical way.

The choice of a Haar system induces an isomorphism between the sheaf of smooth functions on $\mathcal{G}$ and the sheaf of densities along the source fibers, and hence we can transport the convolution product over to the compactly supported functions where it is given by:
\begin{equation*}
(f_1\ast f_2)(g)=\int_{s^{-1}(s(g))}f_1(gh^{-1})f_2(h)\lambda_{s(g)}(h)
\end{equation*}
In particular on the adiabatic groupoid $\Gad$ if we have two compactly supported functions $f_1,f_2$ we obtain:
\begin{equation*}
(f_1\ast f_2)(g,\tau)=|\tau|^{-d}\int_{s^{-1}(s(g))}f_1(gh^{-1},\tau)f_2(h,\tau)\lambda_{s(g)}(h)\,\,\,\,(\tau\neq 0)
\end{equation*}
\begin{equation*}
(f_1\ast f_2)(v,0)=\int_{A_{\pi(v)}}f_1(v-w,0)f_2(w,0)\lambda_{\pi(v)}(w)
\end{equation*}
At this point we notice that the convolution at $\tau=0$ does not require the functions to be compactly supported on $A_x$, being Schwartz is enough (c.f. the usual theory of Fourier transform in $\mathbb{R}^n$). This allows us, in the case of $\Gad$, to enlarge the type of functions/densities on which we let the deformation complex act.

To this end we refer to the work of \cite{cr}, where a Fr\'ech\`et algebra $\mathscr{S}_c(\Gad)$ is constructed with evaluations
\[
\mathscr{S}_c(\Gad)_t=\begin{cases} \mathscr{S}_c(A)& t=0\\ C_c^\infty(\mathcal{G})&t\not = 0.\end{cases}
\]
Here $\mathscr{S}_c(A)$ denotes the space of functions that are Schwartz along the fibers of the Lie algebroid and have compact support along $M$. This Schwartz type algebra should be thought of as a dense subalgebra the reduced $C^*$-algebra $C^*_r(\Gad)$.

By the discussion above, the convolution product is perfectly well-defined on $\mathscr{S}_c(\Gad)$ and we can extend our viewpoint of the map $\chainmap: C^\bullet_\Def(\Gad)\to C^\bullet_\text{Hoch}(\mathcal{A}_{\Gad})$ to let $\chainmap(c)$ (for $c\in C^k_\Def(\Gad)$) act on functions in $\mathscr{S}_c(\Gad)$. At this point it should be remarked that the isomorphism between functions and densities induced by a Haar system does not preserve the action of vector fields (indeed on the level of densities one also needs to compare $\mathcal{L}_X\lambda$ with $\lambda$!). So really we should introduce in parallel to $\mathscr{S}_c(\Gad)$ the notion of densities with are of Schwartz-type along $\tau=0$, but for the sake of not being overly pedantic we will not do this and just be careful when writing down the action of $\chainmap(c)$.

In what follows for a smooth family $\{f_t\}_{t\neq 0}$ of compactly supported functions on $\mathcal{G}$ and $f'\in\mathscr{S}_c(A)$ we will use the notation
\begin{equation*}
\lim_{t\to 0} f_t=f'
\end{equation*}
if the function $F:\Gad\to\mathbb{R}$ given by
\begin{equation*}
F(g,t)=f_t(g)
\end{equation*}
\begin{equation*}
F(v,0)=f'(v)
\end{equation*}
is an element of $\mathscr{S}_c(\Gad)$.
\subsection{Fourier transform on vector bundles}
We briefly discuss the notion of Fourier transform on a vector bundle $E\to M$ under the choice of a Haar system on $E$. This discussion follows the results of Landsman and Ramazan \cite[\S 7]{landsman}. Recall that a vector bundle $\pi: E\to M$ can be seen as a groupoid over $M$ where both the source and the target map are the projection $\pi$ and the multiplication is the fiberwise addition. Since $\ker(d\pi)\cong\pi^\ast E$ a choice of a Haar system is at every $v\in E$ a choice of a density on $E_{\pi(v)}$ that is invariant, where invariance in this case means that the choice is constant along the fiber.

If we choose such a Haar system $\{\mu_x\}_{x\in M}$, in \cite{landsman} the Fourier transform $\mathcal{F}_\mu:\mathscr{S}(E)\to\mathscr{S}(E^\ast)$ was defined by
\begin{equation*}
(\mathcal{F}_\mu f)(\xi_x)=\int_{E_x} f(v)e^{-i\langle\xi_x,v\rangle}d\mu_x(v)
\end{equation*}
Furthermore, it was shown that this map is a linear isomorphism which intertwines the $\mu$-convolution product on $E$ and the pointwise product on $E^\ast$, and when $(x,v)$ are coordinates on $E$ induced by a frame with dual coordinates $(x,\xi)$ we have for $f\in\mathscr{S}(E)$, $g\in\mathscr{S}(E^\ast)$ and $a\in C^\infty(M)$ that
\begin{align*}
\mathcal{F}_\mu((a\circ\pi)f)&=(a\circ \pi)\mathcal{F}_\mu\\
\frac{\partial\mathcal{F}_\mu(f)}{\partial x_j}&=\mathcal{F}_\mu\left(\frac{\partial f}{\partial x_j}\right)+\left(\frac{\partial\text{log}(\mu_e)}{\partial x_j}\circ\pi\right)\mathcal{F}_\mu(f)\\
\frac{\partial\mathcal{F}_\mu(f)}{\partial \xi_j}&=-i\mathcal{F}_\mu(v_j f)\\
\frac{\partial\mathcal{F}_\mu^{-1}(g)}{\partial v_j}&=i\mathcal{F}_\mu^{-1}(\xi_jg)
\end{align*}
Note that after the choice of a Haar system $\mu$ we obtain an isomorphism between the algebra of functions $C^\infty_c(E)$ with the $\mu$-convolution product and the convolution algebra $\mathcal{A}_E$ of densities with the (intrinsic) convolution product. In particular if $X\in\mathfrak{X}(E)$, we can see $\chainmap(X)$ as defined on functions (which is, again, not equal to the usual action of vector fields on functions), and we can extend the action to Schwartz functions.

Now using the Fourier transform, we can transport the action on the convolution algebra of $E$ to an action on the usual algebra with the pointwise product on $E^\ast$.
\begin{prop}
Let $X$ be a linear vector field on $E$. Then the map $\hat{X}:\mathscr{S}(E^\ast)\to\mathscr{S}(E^\ast)$ given by
\begin{equation*}
\hat{X}(f)=\mathcal{F}_\mu(\chainmap(X)(\mathcal{F}_\mu^{-1}(f)))
\end{equation*}
defines a linear vector field on $E^\ast$. Here $\Phi$ is the natural chain map we defined before, applied to the vector bundle $E$ seen as a groupoid.
\end{prop}
\begin{proof}
First we show that $\hat{X}$ is indeed a vector field, i.e. a derivation with respect to the pointwise product. Since $\hat{X}$ is the conjugation of $\chainmap(X)$ with an isomorphism which intertwines the convolution product on $\mathscr{S}(E)$ and the pointwise product on $\mathscr{S}(E^\ast)$ this is equivalent to showing that $\chainmap(X)$ is a derivation for the convolution product. When we see $E\to M$ as a groupoid, this is equivalent to showing that $X$ is a multiplicative vector field, and it is easy to see that on a vector bundle the multiplicative vector fields are precisely the linear vector fields.

To see that $\hat{X}$ is a linear vector field we do a local computation on a trivial vector bundle $E=\mathbb{R}^m_x\times\mathbb{R}^n_v\to\mathbb{R}^m_x$ with Haar system $f(x)dv_1\wedge\cdots\wedge dv_n$. Using the properties of the Fourier transform stated before it follows that if
\begin{equation*}
X(x,v)=\sum_{i=1}^mX_i(x)\frac{\partial}{\partial x_i}+\sum_{j=1}^n\sum_{k=1}^n Y_{jk}(x)v_j\frac{\partial}{\partial v_k}
\end{equation*}
then
\begin{equation*}
\hat{X}(x,\xi)=\sum_{i=1}^mX_i(x)\frac{\partial}{\partial x_i}-\sum_{j=1}^n\sum_{k=1}^nY_{jk}(x)\xi_k\frac{\partial}{\partial\xi_j}
\end{equation*}
which indeed shows that $\hat{X}$ is a linear vector field.
\end{proof}
Recall that a linear vector field $X\in\mathfrak{X}(E)$ is the same as a linear map $X:\Gamma(E^\ast)\to\Gamma(E^\ast)$ with a symbol $s_X\in\mathfrak{X}(M)$ such that
\begin{equation*}
X(f\alpha)=fX(\alpha)+s_X(f)\alpha\qquad(f\in C^\infty(M),~\alpha\in\Gamma(E)).
\end{equation*}
Furthermore, recall the canonical pairing $\langle-,-\rangle:\Gamma(E^\ast)\times\Gamma(E)\to C^\infty(M)$. Then for a linear vector field $X$, the local calculation from the proof above generalizes to the following.
\begin{prop}\label{xhat}
Let $X\in\mathfrak{X}(E)$ be a linear vector field, then the linear vector field $\hat{X}\in\mathfrak{X}(E^\ast)$ is uniquely determined by the fact that for $\beta\in\Gamma(E^\ast)$ and $\alpha\in\Gamma(E)$
\begin{equation*}
\langle\beta,\hat{X}(\alpha)\rangle+\langle X(\beta),\alpha\rangle=s_X(\langle\beta,\alpha\rangle)
\end{equation*}
\end{prop}
We can play a similar game, albeit slightly more involved in notation, for higher order deformation elements of the vector bundle. So consider an element $X\in\tilde{C}^k_\Def(E)$ given by
\begin{equation*}
X(v_1,...,v_n)=X_1(v_1)\langle \beta_2,v_2\rangle\cdots\langle \beta_k,v_k\rangle
\end{equation*}
where $X_1$ is a linear vector field on $E$ and $\beta_2,...,\beta_k\in\Gamma(E^\ast)$. One immediately checks that this is a closed element of $\tilde{C}^k_\Def(E)$, so that the Fourier transform
\begin{equation*}
\hat{X}(f_1,...,f_k)=\mathcal{F}_\mu(\chainmap(X)(\mathcal{F}_\mu^{-1}(f_1),...,\mathcal{F}_\mu^{-1}(f_k)))
\end{equation*}
is a closed element of the Hochschild complex of $C^\infty(E^\ast)$. By the specific form of $X$ is it easy to see that
\begin{equation*}
\chainmap(X)(a_1,...,a_k)=\chainmap(X_1)(a_1)\ast (\beta_2a_2)\ast\cdots\ast(\beta_k a_k)
\end{equation*}
where we see the $s_i$ as fiberwise linear maps on $E$. In particular we see that
\begin{equation*}
\hat{X}=\hat{X_1}\otimes\hat{\beta_2}\otimes\cdots\otimes\hat{\beta_k}
\end{equation*}
where for $\beta\in\Gamma(E^\ast)$, $\hat{\beta}$ is the vector field on $E^\ast$ given by
\begin{equation*}
\hat{\beta}(f)=\mathcal{F}_\mu(\beta\mathcal{F}_\mu^{-1}(f))
\end{equation*}
A local computation shows that $\hat{\beta}$ is identically zero on fiberwise constant maps and for the map induced by a section $\alpha\in\Gamma(E)$ we have
\begin{equation*}
\hat{\beta}(\alpha)=\frac{1}{i}\langle \beta,\alpha\rangle
\end{equation*}
In particular, we see that if we anti-symmetrize, we obtain the linear multivectorfield $\hat{X_1}\wedge\hat{\beta_2}\wedge\cdots\wedge\hat{\beta_k}$ on $E^\ast$.
\subsection{Deformation quantization of $A^\ast$ and the Van Est map}
Now, fix a choice of a Haar system of $\mathcal{G}$, which by the discussion above induces a Haar system on $\Gad$ and a Haar system $\mu$ on $A\to M$. The last one makes sure that we can talk about a Fourier transform $\mathcal{F}_\mu:\mathscr{S}(A)\to\mathscr{S}(A^\ast)$.

Slightly tweaking the results of \cite{landsman} we obtain {\em quantization maps} $q_t:\mathscr{S}_c(A^*)\to C^\infty_c(\mathcal{G}),~t\not = 0$ given by
\[
q_t(f)(g):=\chi(g)\mathcal{F}_\mu^{-1}(f)(\frac{1}{t}\exp^{-1}(g)),
\]
which satisfy
\begin{equation}
\label{quant}
\lim_{t\to 0}(q_t(f_1 f_2)-q_t(f_1)\ast q_t(f_2))=0,\quad \lim_{t\to 0}(\frac{1}{it}[q_t(f_1),q_t(f_2)]-q_t(\{f_1,f_2\}))=0.
\end{equation}
Here $\chi\in C^\infty_c(\mathcal{G})$ is a cut-off function that equals $1$ in a neighborhood of $M\subset\mathcal{G}$ 
with support inside an open neighbourhood of the units onto which the exponential map is a diffeomorphism. The Poisson 
bracket $\{~,~\}$ is the bracket associated to the so-called Lie--Poisson structure on $A^*$.

Explicitely, one of the differences with the results of \cite{landsman} is that we do not need the property $q_t(f^\ast)=q_t(f)^\ast$ for which the Weyl exponential map $\text{exp}^{\text{W}}$ is used, and instead we can use the normal exponential map. Secondly, we do not need to restrict to Paley-Wiener functions, as we allow for Schwarz-type functions at $t=0$ and use the cut-off function on the level of $\mathcal{G}$ instead of $A$, the deviation vanishing as $t$ approaches $0$. Lastly, as the relevant calculations on the local forms in $A$ and $A^\ast$ are valid for all Schwarz functions and not just Paley-Wiener functions, the relevant propositions in \cite{landsman} still hold in this situation. The variety of quantizations by using different types of exponential maps is also reflected on the more algebraic level in \cite{nw} by using different orderings in the Fedossov construction of {\em formal} deformation quantizations of $A^*$.

We now briefly recall the van Est-map as given in \cite[\S 10]{cms}. First the deformation complex of the algebroid $C^k_\Def(A)$ is given by antisymmetric multilinear maps $D: \Gamma(A)^k\to\Gamma(A)$ that have a symbol $s_D:\Gamma(A)^{k-1}\to\mathfrak{X}(M)$ such that
\begin{equation*}
D(\alpha_1,...,f\alpha_k)=fD(\alpha_1,...,\alpha_k)+s_D(\alpha_1,...,\alpha_{k-1})(f)\alpha_k
\end{equation*}
Note that we can, and will, see elements of $C^k_\Def(A)$ as linear multivectorfields on $A^\ast$ (by noting that sections of $A$ are the same as fiberwise linear maps on $A^\ast$) and in turn see the deformation complex of $A$ as the linear Poisson complex of the Poisson manifold $A^\ast$.

Then for $\alpha\in\Gamma(A)$ there are maps $R_\alpha:\tilde{C}^k_\Def(\mathcal{G})\to\tilde{C}^{k-1}_\Def(\mathcal{G})$ which are given for $k=1$ by
\begin{equation*}
R_\alpha(c)=[c,\overrightarrow{\alpha}]|_M
\end{equation*}
and for $k>0$ by
\begin{equation*}
R_\alpha(c)(g_1,...,g_{k-1})=(-1)^{k-1} \left.\frac{d}{d\epsilon}\right|_{\epsilon=0} c(g_1,...,g_{k-1},\Phi^\epsilon_{\overrightarrow{\alpha}}(s(g_{k-1}))^{-1})
\end{equation*}
The van Est-map $\mathcal{V}:\tilde{C}^k_\Def(\mathcal{G})\to C^k_\Def(A)$ is then given by
\begin{equation*}
\mathcal{V}(c)(\alpha_1,...,\alpha_k)=\sum_{\sigma\in S_k}(-1)^\sigma (R_{\alpha_{\sigma(k)}}\circ\cdots\circ R_{\alpha_{\sigma(1)}})(c)
\end{equation*}
The connection between the van Est-map and the quantization maps is then as follows.
\begin{thm}
Let $k\geq 1$ and $c\in \tilde{C}^k_{\Def}(\mathcal{G})$ and suppose the choice of a Haar system on $\mathcal{G}$ inducing a Haar system $\mu$ on the algebroid $A$. Given $f_1,\ldots,f_k\in  \mathscr{S}_c(A^*)$, the following 
equality holds true:
\[
\mathcal{V}(c)(f_1,\ldots,f_k)=\mathcal{F}_\mu\left(\lim_{t\to 0}\left(\sum_{\sigma\in S_k}(-1)^\sigma\frac{1}{(it)^{k-1}}\chainmap(c)(q_t(f_{\sigma(1)}),\ldots,q_t(f_{\sigma(k)}))\right)\right)
\]
\end{thm}
\begin{rmk}
Note that the right hand side of the equation above is well-defined, since the sum of which we take the limit is a function on the groupoid. The limit is then a Schwartz function on the algebroid and so if we take the Fourier transform we obtain a Schwartz function on the dual of the algebroid.
\end{rmk}
\begin{proof}
We start with the case $k=1$. First note that for $f\in\mathscr{S}_c(A^\ast)$ the map $q(f):\Gad\to\mathbb{R}$ given by
\begin{align*}
q(f)(g,t)&=q_t(f)(g)\\
q(f)(v,0)&=\mathcal{F}_\mu^{-1}(f)(v)
\end{align*}
is an element of $\mathscr{S}_c(\Gad)$. Then note that the family $\{c\}_{t\neq 0}$ is a family of vector fields on $\mathcal{G}$ which can be extended to a vector field on $\Gad$, namely to the vector field $c_\text{inv}$ obtaines by \cref{vfonnms}. Then notice that $\chainmap(c_\text{inv})(q(f))$ is an element of $\mathscr{S}_c(\Gad)$ consisting of
\begin{align*}
\chainmap(c_\text{inv})(q(f))_t&=\chainmap(c)(q_t(f)),\qquad(t\neq 0)\\
\chainmap(c_\text{inv})(q(f))_0&=\chainmap(c_0)(\mathcal{F}_\mu^{-1}(f))
\end{align*}
where $c_0$ is the linear vector field that is the restriction of $c_\text{eqv}$ to $t=0$. Note that it is linear, since it is the application of \cref{vectorfieldnormalbundle} to the vector field $c$ on $\mathcal{G}$. In particular we see that
\begin{equation*}
\lim_{t\to 0}\chainmap(c)(q_t(f))=\chainmap(c_0)(\mathcal{F}_\mu^{-1}(f))
\end{equation*}
and so we need to show that $\mathcal{V}(c)=\hat{c_0}$.

By \cref{xhat} this means that we need to show that for $\beta\in\Gamma(A^\ast)$ and $\alpha\in\Gamma(A)$ we have
\begin{equation*}
\langle\beta,\mathcal{V}(c)(\alpha)\rangle+\langle c_0(\beta),\alpha\rangle=s_c(\langle\beta,\alpha\rangle)
\end{equation*}
We note two things. First that every $\beta\in\Gamma(A^\ast)=C^\infty_\text{lin}(A)$ can be written as $d_nh$ for $h\in C^\infty(\mathcal{G})$ with $h|_M=0$. Second that since we have an explicit inclusion of $A$ into the tangent bundle of $\mathcal{G}$ this means that:
\begin{equation*}
\langle d_nh,\alpha\rangle(x)=\alpha(x)(h)
\end{equation*}
We are now ready to show the equality. First we have
\begin{align*}
\langle d_nh,\mathcal{V}(c)(\alpha)\rangle(x)&=[c,\overrightarrow{\alpha}](1_x)(h)\\
&=c(1_x)(\overrightarrow{\alpha}(h))-\alpha(x)(c(h))
\end{align*}
\begin{equation*}
\langle c_0(d_n h),\alpha\rangle(x)=\langle d_n(ch),\alpha\rangle(x)=\alpha(x)(c(h))
\end{equation*}
and since $c_{1_x}=du(s_c(x))$ combined with $\overrightarrow{\alpha}|_M=\alpha$ we have
\begin{equation*}
s_c(\langle d_nh,\alpha\rangle)(x)=s_c(\overrightarrow{\alpha}(h)|_M)(x)=c(1_x)(\overrightarrow{\alpha}(h))
\end{equation*}
For $k>1$ we restrict to the case where $c=c_1\otimes h_2\otimes\cdots\otimes h_k$ with $c_1\in\mathfrak{X}(\mathcal{G})$ and $h_2,...,h_k\in C^\infty(\mathcal{G})$. For $c$ to be an element of $\tilde{C}^k_\Def(\mathcal{G})$ it is necessairy and sufficient to have $c_1\in\tilde{C}^1_\Def(\mathcal{G})$ and $h_i|_M=0$. Similar to the case $k=1$ we note that
\begin{equation*}
\frac{1}{(it)^{k-1}}\chainmap(c)(q_t(f_1),...,q_t(f_k))=\chainmap(\frac{1}{(it)^{k-1}}c)(q_t(f_1),...,q_t(f_k))
\end{equation*}
which, as $t\to 0$, converges
\begin{equation*}
\chainmap(c_0)(\mathcal{F}_\mu^{-1}(f_1),...,\mathcal{F}_\mu^{-1}(f_k))
\end{equation*}
if we find a vector field $c_0$ on $A$ that together with the family $\{\frac{1}{(it)^{k-1}}c\}_{t\neq 0}$ defines a smooth deformation element of $\Gad$.

To calculate this localization we remark that we can do the calculation in $\mathcal{G}^k$ using the cartesian product of the exponential map $A\to\mathcal{G}$, in stead of working in $\mathcal{G}^{(k)}$ and using the machinery of the previous section. This is for two reasons: firstly our definition of $c$ extends to $\mathcal{G}^k$. Secondly the difference of $(v_1,...,v_k)\in A^{\oplus k}$ seen as tangent vectors on $\mathcal{G}^k$ and $(v_1,...,v_k)\in A^{\oplus k}$ seen as tangent vectors in $\mathcal{G}^{(k)}$ which are normal to the units, using the isomorphism of \cref{normalbundenerve}, are tangent vectors in $\mathcal{G}^k$ which are along the units. Since $c$ vanishes along the units, we can neglect this.

Now to do the actual calculation we consider the chart $\theta: A^{\oplus k}\times\mathbb{R}^\ast\to\mathcal{G}^k\times\mathbb{R}^\ast$ given by
\begin{equation*}
\theta(v_1,...,v_k,t)=(\exp(tv_1),...,\exp(tv_k),t)
\end{equation*}
Then if we look at the family $\{\frac{1}{(it)^{k-1}}c\}_{t\neq 0}$, we see that if we take the pullback along $\theta$ we obtain:
\begin{equation*}
\theta^\ast(\{\frac{1}{(it)^{k-1}}c\}_{t\neq 0})(v_1,...,v_k,t)=\frac{1}{(it)^k}c_1(\exp(tv_1))h_2(\exp(tv_2))\cdots h_k(\exp(tv_k))
\end{equation*}
Distributing the $k$ powers of $\frac{1}{t}$ over the $k$ different terms we see that
\begin{equation*}
c_0(v_1,...,v_k)=\frac{1}{i^{k-1}}(c_1)_0(v_1)d_nh_2(v_2)\cdots d_nh_2(v_k)
\end{equation*}
since
\begin{equation*}
\frac{1}{t}c_1(\exp(tv_1))\to (c_1)_0(v_1)
\end{equation*}
\begin{equation*}
\frac{1}{t}h(\exp(tv))\to d_nh(v)
\end{equation*}
as $t\to 0$, so we see that $c_0=\frac{1}{i^{k-1}}(c_1)_0\otimes d_nh_2\otimes\cdots\otimes d_nh_k$, which is a linear deformation element, and we want to show that $\mathcal{V}(c)$ is the anti-symmetrization of the Fourier transform $\hat{c_0}$. By the discussion at the end of the previous subsection we see that $\hat{c_0}$ is determined for $\alpha_1,...,\alpha_k\in\Gamma(A)$ by
\begin{equation*}
\hat{c_0}(\alpha_1,...,\alpha_k)=\frac{1}{i^{2(k-1)}}\hat{(c_1)_0}(\alpha_1)\langle d_nh_2,\alpha_2\rangle\cdots\langle d_nh_k,\alpha_k\rangle
\end{equation*}
Next we investigate $R_\alpha(c)$, we obtain:
\begin{align*}
R_\alpha(c)(g_1,...,g_{k-1})&=(-1)^{k-1}\frac{d}{d\epsilon}|_{\epsilon=0}c_1(g_1)h_2(g_2)\cdots h_{k-1}(g_{k-1})h_k(\Phi^\epsilon_{\overrightarrow{\alpha}}(s(g_k))^{-1})\\
&=(-1)^{k-1}c_1(g_1)h_2(g_2)\cdots h_{k-1}(g_{k-1})dh_k(d\iota(\alpha(s(g_k)))
\end{align*}
Then since $f_k|_M=0$ and for $v\in A_x$ we have $d\iota v=-v+d(u\circ t)(v)$ we obtain
\begin{equation*}
R_\alpha(c)(g_1,...,g_{k-1})=(-1)^k c_1(g_1)h_2(g_2)\cdots h_{k-1}(g_{k-1})d_nh_k(\alpha(s(g_{k-1}))
\end{equation*}
Doing this inductively, and using that the flow of $\overrightarrow{\alpha}$ preserves source fibers, we see
\begin{equation*}
(R_{\alpha_2}\circ\cdots\circ R_{\alpha_k})(c)(g)=(-1)^{\frac{(k-1)(k-2)}{2}}c_1(g)d_nh_2(\alpha_2(s(g))\cdots d_nh_k(\alpha_k(s(g))
\end{equation*}
Then since this is simply $c_1$ multiplied with a function that is constant along the $s$-fibers, we obtain:
\begin{align*}
(R_{\alpha_1}\circ\cdots \circ R_{\alpha_k})(c)&=(-1)^{\frac{(k-1)(k-2)}{2}}\mathcal{V}(c_1)(\alpha_1)\langle d_nh_2,\alpha_2\rangle\cdots\langle d_nh_k,\alpha_k\rangle\\
&=i^{(k-1)(k-2)}\mathcal{V}(c_1)(\alpha_1)\langle d_nh_2,\alpha_2\rangle\cdots\langle d_nh_k,\alpha_k\rangle
\end{align*}
Since already know by the calculation in the case $k=1$ that $\mathcal{V}(c_1)(\alpha_1)=\hat{(c_1)_0}(\alpha_1)$ we see that
\begin{equation*}
(R_{\alpha_1}\circ\cdots \circ R_{\alpha_k})(c)=i^{k(k-1)}\hat{c_0}(\alpha_1,...,\alpha_k)
\end{equation*}
Then note that there is a mismatch in the summation over $S_k$ in $\mathcal{V}(c)$ and in the right hand side of the theorem. In particular the right hand side in the last equation corresponds to the identity permutation in the statement of the theorem, while the right hand side corresponds to the permutation in the definition of $\mathcal{V}(c)$ that sends $j$ to $k-j$. This sign of this permutation is $(-1)^{\frac{k(k-1)}{2}}$, for which we have to correct, so that we obtain
\begin{align*}
\mathcal{V}(c)(\alpha_1,...,\alpha_k)&=\sum_{\sigma\in S_k}(-1)^\sigma (R_{\alpha_{\sigma(k)}}\circ\cdots R_{\alpha_{\sigma(1)}})(c)\\
&=\sum_{\sigma\in S_k}(-1)^\sigma i^{k(k-1)}(R_{\alpha_{\sigma(1)}}\circ\cdots\circ R_{\alpha_{\sigma(k)}})(c)\\
&=\sum_{\sigma\in S_k}(-1)^\sigma i^{2k(k-1)}\hat{c_0}(\alpha_{\sigma(1)},...,\alpha_{\sigma(k)})\\
&=\sum_{\sigma\in S_k}(-1)^\sigma \hat{c_0}(\alpha_{\sigma(1)},...,\alpha_{\sigma(k)})
\end{align*}
So we see that $\mathcal{V}(c)$ equals the linear multivector field that is the antisymmetrization of $\hat{c_0}$. In particular this means that for $f_1,...,f_k\in\mathscr{S}_c(A^\ast)$ we have
\begin{align*}
\mathcal{V}(c)(f_1,...,f_k)&=\sum_{\sigma\in S_k}(-1)^\sigma\hat{c_0}(f_{\sigma(1)},...,f_{\sigma(k)})\\
&=\frac{1}{i^{k-1}}\mathcal{F}_\mu\left(\lim_{t\to 0}\left(\sum_{\sigma\in S_k}(-1)^\sigma\frac{1}{(it)^{k-1}}\chainmap(c)(q_t(f_{\sigma(1)}),\ldots,q_t(f_{\sigma(k)}))\right)\right)
\end{align*}
This completes the proof.
\end{proof}

\begin{rmk}
This theorem, restricted to multiplicative vector fields, can be viewed as a statement about the ``classical limit'' of certain derivations of the convolution algebra, and looks very similar to certain aspects of the proof of the Atiyah--Singer index theorem given in \cite{enn}. Indeed, it would be interesting to investigate its use in index theory for Lie groupoids, as it 
exactly fits into the framework of relating the van Est map to the classical limit, as shown in the index theorem of \cite{ppt} 
for smooth groupoid cohomology $H^\bullet_{\rm diff}(\mathcal{G})$.
\end{rmk}

In the previous proof we have only used the fact that $q_t(f)$ converges to $\mathcal{F}_\mu^{-1}(f)$ in $\mathscr{S}_c(\Gad)$ as $t$ goes to $0$, we have not used the properties which makes the family $\{q_t\}_{t\neq 0}$ a family of quantization maps, namely their compatibility with the Poisson bracket. However, we have not introduced these specific maps without reason, since we will use the fact that
\begin{equation*}
\lim_{t\to 0}(\frac{1}{it}[q_t(f_1),q_t(f_2)])=\lim_{t\to 0}q_t(\{f_1,f_2\}))
\end{equation*}
to give an alternative proof of the fact that the Van Est map is a {\em chain map}, i.e, compatible with the differentials:
\begin{cor}
The van Est map $\mathcal{V}:\tilde{C}_\Def^\bullet(\mathcal{G})\to C^\bullet_{\text{Pois,lin}}(A^\ast)$ is a chain map.
\end{cor}
\begin{proof}
Let $c\in\tilde{C}^k_\Def(\mathcal{G})$ for $k\geq 1$ and we start by dissecting $\mathcal{V}(\partial c)$. Using the previous theorem we obtain
\small
\begin{align*}
\mathcal{V}(\delta c)(f_1,...,f_{k+1})=&\mathcal{F}_\mu\left(\lim_{t\to 0}\left(\sum_{\sigma\in S_{k+1}}(-1)^\sigma\frac{1}{(it)^{k}}\chainmap(\delta c)(q_t(f_{\sigma(1)}),\ldots,q_t(f_{\sigma(k+1)}))\right)\right)\\
=&\mathcal{F}_\mu\left(\lim_{t\to 0}\left(\sum_{\sigma\in S_{k+1}}(-1)^\sigma\frac{1}{(it)^{k}}(\delta_\text{Hoch}\chainmap(c))(q_t(f_{\sigma(1)}),\ldots,q_t(f_{\sigma(k+1)}))\right)\right)\\
=&\mathcal{F}_\mu\left(\lim_{t\to 0}\left(\sum_{\sigma\in S_{k+1}}(-1)^\sigma\frac{1}{(it)^k}[q_t(f_{\sigma(1)}),\chainmap(c)(q_t(f_{\sigma(2)}),...,q_t(f_{\sigma(k+1)}))]\right)\right)\\
&+\mathcal{F}_\mu\left(\lim_{t\to 0}\left(\sum_{j=1}^k\sum_{\substack{\sigma\in S_{k+1}\\\sigma^{-1}(j)<\sigma^{-1}(j+1)}}(-1)^\sigma (-1)^j\frac{1}{(it)^k}\chainmap(c)(q_t(f_{\sigma(1)}),...,[q_t(f_{\sigma(j)}),q_t(f_{\sigma(j+1)})],...,q_t(f_{\sigma(k)}))\right)\right)
\end{align*}
\normalsize
Now the relation between the commutator, the Poisson bracket and the quantization maps, we can use 1 power of $\frac{1}{it} $ to turn the commutators into Poisson brackets. Also using the fact that $q_t(f)\to \mathcal{F}_\mu^{-1}(f)$ as $t\to 0$ this results in
\small
\begin{align*}
\mathcal{V}(\delta c)(f_1,...,f_{k+1})=&\sum_{\sigma\in S_{k+1}}(-1)^\sigma\left\{f_{\sigma(1)},\mathcal{F}_\mu\left(\lim_{t\to 0}\left(\frac{1}{(it)^{k-1}}\chainmap(c)(q_t(f_{\sigma(2)}),...,q_t(f_{\sigma(k+1)}))\right)\right)\right\}\\
&+\mathcal{F}_\mu\left(\lim_{t\to 0}\left(\sum_{j=1}^k\sum_{\substack{\sigma\in S_{k+1}\\\sigma^{-1}(j)<\sigma^{-1}(j+1)}}(-1)^\sigma (-1)^j\frac{1}{(it)^{k-1}}\chainmap(c)(q_t(f_{\sigma(1)}),...,q_t(\{f_{\sigma(j)},f_{\sigma(j+1)}\}),...,q_t(f_{\sigma(k)}))\right)\right)
\end{align*}
\normalsize
Then using the previous Theorem in reverse order we see that this leads to
\begin{align*}
\mathcal{V}(\delta c)(f_1,...,f_{k+1})&=\sum_{j=1}^{k+1}(-1)^{j+1}\left\{f_j,\mathcal{V}(c)(f_1,...,\hat{f_j},...,f_{k+1})\right\}\\
&+\sum_{j_1<j_2}(-1)^{j_1+j_2}\mathcal{V}(c)(\{f_{j_1},f_{j_2}\},f_1,...,\hat{f_{j_1}},\hat{f_{j_2}},...,f_{k+1})
\end{align*}
which shows that the Van Est map is a chain map.
\end{proof}

\end{document}